\newif\ifarxived
\arxivedtrue
\newif\ifextended
\ifx\Extended\Kaba\else\extendedtrue\fi
\newif\ifprivate
\ifx\Private\Kaba\else\privatetrue\extendedtrue\fi
\newif\ifJapanese
\newif\iftesting
\ifextended\testingtrue\fi
\ifarxived
\extendedtrue
\testingfalse
\privatefalse
\fi
\documentclass[12pt]{article}
\ifextended
\ifarxived
\usepackage[bookmarks=true,backref=page,bookmarksnumbered=true, 
  bookmarkstype=toc,
  pdftitle={Weakly extendible cardinals and compactness of extended logics},%
  pdfauthor={Sakae Fuchino, Hiroshi Sakai}
  pdfsubject={},%
  pdfkeywords={Loewenheim-Skolem theorems, compactness, extendible cardinals, subtle 
    cardinals, strongly uplifting cardinals, weakly extendible 
    cardinals, strongly unfoldable cardinals, shewd cardinals, weakly 
    compact cardinals}]{hyperref}
\else
\usepackage[dvipdfmx,backref=page,bookmarks=true,bookmarksnumbered=true,
  bookmarkstype=toc,
  pdftitle={Weakly extendible cardinals and compactness of extended logics},%
  pdfauthor={Sakae Fuchino (), Hiroshi Sakai ()},%
  pdfsubject={},%
  pdfkeywords={compactness, extendible cardinals, shewd cardinals, strongly unfoldable cardinals, weakly 
    compact cardinals}]{hyperref}
\usepackage[dvipdfmx]{pxjahyper}
\hypersetup{
  backref=page,
  colorlinks=true,
  linkcolor=blue,
  filecolor=blue,      
  citecolor=cyan,
  urlcolor=cyan
}
\fi 
\else 
\newcommand{\href}[2]{#2}
\def\hyperref[#1]#2{#2}
\newcommand{\phantomsection}{}
\fi

\ifarxived
\usepackage{graphicx, color}
\else
\usepackage[dvipdfmx]{graphicx, color}
\fi
\usepackage{amsmath, amssymb}
\usepackage{bbm}
\usepackage{dsfont}
\usepackage{bbold}
\usepackage{marginnote}
\usepackage{mathtools}
\usepackage[normalem]{ulem}
\usepackage{setspace}

\definecolor{darkelectricblue}{rgb}{0.33, 0.39, 0.52}
\definecolor{darkgreen}{rgb}{0.31, 0.47, 0.26}
\newcommand{\extendedcolor}{\color{darkelectricblue}}
\newcommand{\darkgreen}{\color{darkgreen}}
\newcommand{\privatecolor}{\color{darkgreen}}
\ifextended
\newcommand{\darkred}{\color[rgb]{0.8,0.1,0.1}}
\newcommand{\It}{\it\darkred{}}
\else
\newcommand{\darkred}{}
\newcommand{\It}{\it{}}
\fi
\usepackage{theorem}
\newcommand{\bbd}[1]{{\mathbb{#1}}}
\theorembodyfont{\ifJapanese\rm\else\it\fi}
\newcount\minute	
\newcount\hour		
\newcount\hourMins  
\ifJapanese
\def\today%
{
  \the\year "N\,\zeroPadTwo{\the\month}ŒŽ\,\zeroPadTwo{\the\day}"ú%
}
\fi
\def\now%
{
  \minute=\time    
  \hour=\time \divide \hour by 60 
  \hourMins=\hour \multiply\hourMins by 60
  \advance\minute by -\hourMins 
  \zeroPadTwo{\the\hour}:\zeroPadTwo{\the\minute}%
}
\def\zeroPadTwo#1%
{
  \ifnum #1<10 0\fi    
  #1
}
\title{Weakly extendible cardinals and compactness of extended logics}
\ifarxived
\author{\protect\scalebox{1}[1.4]{\frame{\rule{1em}{0ex}\rule{0em}{1.8ex}}\,%
    \frame{\rule{1em}{0ex}\rule{0em}{1.8ex}}\ \ 
    \frame{\rule{1em}{0ex}\rule{0em}{1.8ex}}\qquad\qquad\qquad\quad 
    \frame{\rule{1em}{0ex}\rule{0em}{1.8ex}}\,\frame{\rule{1em}{0ex}\rule{0em}{1.8ex}}\ \frame{\rule{1em}{0ex}\rule{0em}{1.8ex}}\,\frame{\rule{1em}{0ex}\rule{0em}{1.8ex}}}\medskip\\
  \protect\scalebox{1}[1.4]{\quad Saka\'e Fuchino$^{*,\dagger}$\quad and \quad\ Hiroshi Sakai$^{*,\ddagger}$}\\
  }
\else
\author{\protect\scalebox{1}[1.4]{\ \qquad\qquad\qquad\quad }\medskip\\
  \protect\scalebox{1}[1.4]{\quad Saka\'e Fuchino$^{*,\dagger}$\quad and \quad\ Hiroshi Sakai$^{*,\ddagger}$}\\
  }
\fi
\date{}

\ifprivate\fi
\renewcommand{\baselinestretch}{1.2}
\renewcommand{\thefootnote}{(\arabic{footnote})\,}
\iftesting
\fi
\iftesting
\newcommand{\Label}[1]{\label{#1}\marginnote{{\color{cyan}\renewcommand{\baselinestretch}{0.4}\tiny 
		  \rlap{#1}}}}

\else
\newcommand{\Label}[1]{\label{#1}}

\fi

\ifarxived
\def\memo#1{\ifprivate\marginnote{{\darkgreen\normalsize%
      \color{cyan}\mbox{}\vspace{-2.52ex}\\
      \mbox{}\rlap{\scalebox{0.72}[1]{\parbox{15em}{\tiny #1}}}}}
			\fi }%
\else
\def\memo#1{\ifprivate\marginnote{{\darkgreen\normalsize%
      \color{cyan}\mbox{}\vspace{-2.52ex}\\
      \mbox{}\rlap{\scalebox{0.72}[1]{\parbox{15zw}{\tiny #1}}}}}
			\fi }%
\fi
\def\oldmemo#1{}
\def\imemo#1{\ifprivate{\color{cyan}\footnotesize\renewcommand{\baselinestretch}{0.4}\tiny\mbox{}\\%
			#1\par\ignorespaces }\fi }%
\def\oldimemo#1{}
\newcounter{frml}[section]
\newcounter{frmla}[section]
\def\thefrml{{\arabic{section}.\arabic{frml}}}
\def\thefrmla{{$\aleph$\arabic{section}.\arabic{frmla}}}
\def\frmlabel#1{\refstepcounter{frml}{\def\baka{#1}\ifx\baka\empty\else\label{#1}\fi}%
{\rm({\thefrml})\hfill\hfill\hfill}}
\def\ifrmlabel#1{\refstepcounter{frml}{\def\baka{#1}\ifx\baka\empty\else\label{#1}\fi}%
{\rm({\thefrml})\,:\ifarxived\hspace{0.4em}\else\hspace{0.4zw}\fi}}
\def\frmlabela#1{\refstepcounter{frmla}{\def\baka{#1}\ifx\baka\empty\else\label{#1}\fi}%
{\rm({\thefrmla})\hfill\hfill\hfill}}
\def\xitem[#1]{\item[\frmlabel{#1}]\mbox{}%
	\iftesting\marginnote{\color{cyan}{\mbox{\tiny#1}}}\fi\ignorespaces}
\def\xitemq[#1]{\item[\frmlabel{#1}]\mbox{}%
	\ignorespaces}
\def\xitemd[#1]#2{\item[(\ref{#1})$#2$\hfill\hfill\hfill]}
\iftesting
\def\ixitem[#1]{\ifrmlabel{#1}\marginnote{{\color{cyan}\renewcommand{%
				\baselinestretch}{0.6}\hfill\tiny \rlap{#1}\mbox{}}}\ignorespaces}
\def\ixitemx[#1]{\ifrmlabel{#1}\marginnote{\mbox{}\vspace{-1zh}\\{\color{cyan}\renewcommand{%
				\baselinestretch}{0.6}\hfill\tiny \rlap{#1}\mbox{}}}\ignorespaces}
\else
\def\ixitem[#1]{\ifrmlabel{#1}\ignorespaces}
\def\ixitemx[#1]{\ifrmlabel{#1}\ignorespaces}
\fi
\def\xitemA[#1]{\item[\frmlabela{#1}]\mbox{}%
	\iftesting\marginnote{\color{cyan}\mbox{\tiny#1}}\fi\ignorespaces}
\def\xitemsub[#1]#2{\item[\frmlabel{#1}$_{#2}$]\mbox{}%
	\iftesting\marginnote{\color{cyan}{\renewcommand{%
				\ba
selinestretch}{0.6}\tiny#1}}\fi\ignorespaces}
\def\xxitem[#1][#2]{\item[(\ref{#1}{\makebox[1.4ex][c]{#2}})]\mbox{}%
	\iftesting\marginnote{\color{cyan}{\renewcommand{%
				\baselinestretch}{0.6}\tiny\{#1\}\{#2\}}}\fi\ignorespaces}
\def\xitemof#1{{\rm({\ref{#1}})}}
\def\Xitem[#1]{\item[{\makebox[7ex][l]{\rm(\ref{#1})}}]\iftesting\marginnote{{\color{cyan}\renewcommand{%
				\baselinestretch}{0.6}\tiny#1}}\fi\ignorespaces}

\def\xitemdof#1{$\mbox{\rm(\ref{#1})}'$}

\newenvironment{xitemize}{\begin{list}{}{\parsep=0.5\smallskipamount%
			\itemindent=-0.4ex%
			\itemsep=0.5\smallskipamount\leftmargin=4em\labelwidth=3em\labelsep=0.7em}}%
							 {\end{list}}
\def\assert#1{\noindent\makebox[4.8ex][r]{\rm(\makebox[2.2ex][c]{#1})}\ \ \ignorespaces}
\def\wassert#1{\assert{#1}}
\def\wassertof#1{\makebox[4.8ex][r]{\rm(\makebox[2.2ex][c]{#1})}}%
\def\assertof#1{(#1)}%

\def\daimaru#1{\makebox[1em][c]{\mbox{\leavevmode\lower.144ex\hbox{
        \rlap{\hbox to 
          0.76em{\hfil\mbox{}\hfill{}\raisebox{0.054ex}{\ifarxived\scalebox{1.2}{$\bigcirc$}\else\scalebox{1.2}{›}\fi}\hfil}}
        \raise0.342ex\hbox to 1em{\hfil{\hspace{0.16em}\footnotesize#1}\hfil}}}}\,}

\newtheorem{Thm}{\ifJapanese{\bf '藝}\else {\bf Theorem}\fi}[section]

\newtheorem{ThmA}{\ifJapanese{\bf '藝\,A\!}\else{\bf Theorem\,A\!}\fi}[section]

\newtheorem{Prop}[Thm]{\ifJapanese{\bf –½'è}\else{\bf Proposition}\fi}
\newtheorem{PropA}[ThmA]{\ifJapanese{\bf –½'è\,A\!}\else{\bf Proposition\,A\!}\fi}
\newtheorem{Problem}[Thm]{\ifJapanese{\bf –¢‰ðŒˆ–â'è}\else{\bf Problem}\fi}

\newtheorem{Lemma}[Thm]{\ifJapanese{\bf •â'è}\else{\bf Lemma}\fi}
\newtheorem{LemmaA}[ThmA]{\ifJapanese{\bf •â'è\,A\!}\else{\bf Lemma\,A\!}\fi}
\newtheorem{Cor}[Thm]{\ifJapanese{\bf Œn}\else{\bf Corollary}\fi}

\newtheorem{Claim}{{\bf Claim}}[Thm]

\newcommand{\footnotemarkA}{{\addtocounter{footnote}{-1}\extendedcolor\footnotemark${}^{\!\mathrm a}$}}
\newcommand{\footnotetextA}[1]{\footnotetext{\extendedcolor{}${}^{\!\mathrm a}$\, #1}}%

\newcommand{\prf}{\ifJapanese{\bf Ø–¾D\ }\ignorespaces\else{\bf 
		Proof.\ \ }\ignorespaces\fi}

\newcommand{\prfofClaim}{\raisebox{-.4ex}{\Large $\vdash$\ \ }}

\newcommand{\Thmof}[1]{\ifJapanese{'藝\,\ref{#1}}\else{Theorem~\ref{#1}}\fi}

\newcommand{\bfThmof}[1]{\ifJapanese{\bf '藝\,\ref{#1}}\else{\bf Theorem~\ref{#1}}\fi}
\newcommand{\Lemmaof}[1]{\ifJapanese{•â'è\,\ref{#1}}\else{Lemma~\ref{#1}}\fi}

\newcommand{\LemmaAof}[1]{\ifJapanese{•â'è\,A\,\ref{#1}}\else{Lemma\,A\,\ref{#1}}\fi}

\newcommand{\Propof}[1]{\ifJapanese{–½'è\,\ref{#1}}\else{Proposition~\ref{#1}}\fi}
\newcommand{\PropAof}[1]{\ifJapanese{–½'è\,A\,\ref{#1}}\else{Proposition\,A\,\ref{#1}}\fi}

\newcommand{\Claimof}[1]{{Claim \ref{#1}}}

\newcommand{\sectionof}[1]{\ifJapanese{'æ\ref{#1}ß}\else{Section~\ref{#1}}\fi}

\newcommand{\footnoteof}[1]{\ifJapanese{‹r'\ref{#1}}\else{footnote~\ref{#1}}\fi}

\newcommand{\Lemmaabove}{{\ifJapanese •â'è\else Lemma\fi\ \number\theThm}}

\newcommand{\ubecause}[3]{\underbrace{{}#1{}%
  \ifx\bakakaba#2\bakakaba\rule[-0.72ex]{0pt}{1pt}\else\rule[#2]{0pt}{1pt}\fi}_{\footnotesize\clap{#3}}}
\newcommand{\obecause}[3]{\overbrace{{}#1{}%
  \ifx\bakakaba#2\bakakaba\rule[1.62ex]{0pt}{1pt}\else\rule[#2]{0pt}{1pt}\fi}^{\footnotesize\clap{#3}}}
\newsavebox{\qedbox}\sbox{\qedbox}{
{\unitlength=0.05mm \begin{picture}(40,60)
\put(0,0){\framebox(30,44)[cc]{}}
\put(30,-7){\rule{7\unitlength}{44\unitlength}}
\put(10,-7){\rule{27\unitlength}{7\unitlength}}
\end{picture}}}
\newcommand{\qed}{\mbox{}\hfill\usebox{\qedbox}}
\newcommand{\smallqed}%
{\mbox{}\smallskip\hfill\raisebox{-.4ex}{\Large $\dashv$}}
\newcommand{\qedof}[1]%
{\mbox{} \hspace*{\fill}{\usebox{\qedbox}{\tiny~(#1)}}}
\newcommand{\Qedof}[1]%
{\mbox{} \hspace*{\fill}{\usebox{\qedbox}%
{\tiny~(#1~\number\theThm)}}}
\newcommand{\QedAof}[1]%
{\mbox{} \hspace*{\fill}{\usebox{\qedbox}%
{\tiny~(#1~\number\theThmA)}}}
\newcommand{\qedofThm}{\Qedof{\ifJapanese '藝\else Theorem\fi}}

\newcommand{\qedofCor}{\Qedof{\ifJapanese Œn\else Corollary\fi}}
\newcommand{\qedofProp}{\Qedof{\ifJapanese –½'è\else Proposition\fi}}
\newcommand{\qedofLemma}{\Qedof{\ifJapanese •â'è\else Lemma\fi}}

\newcommand{\qedskip}{\medskip}

\newcommand{\qedofClaim}%
{\mbox{}\hfill\raisebox{-.4ex}{\Large $\dashv$ }\nolinebreak%
\mbox{\tiny~(Claim~\number\theClaim)}}
\newcommand{\qedofClaimA}%
{\mbox{}\hfill\raisebox{-.4ex}{\Large $\dashv$ }\nolinebreak%
\mbox{\tiny~(Claim~A\,\number\theClaimA)}}
\newcommand{\qedofClaimAof}[1]%
{\mbox{}\hfill\raisebox{-.4ex}{\Large $\dashv$ }\nolinebreak%
\mbox{\tiny~(Claim~A\,\ref{#1})}}
\newcommand{\qedofSubclaim}%
{\mbox{}\hfill\raisebox{-.4ex}{\Large $\dashv$ }\nolinebreak%
\mbox{\tiny~(Subclaim~\number\theSubclaim)}}
\newcommand{\cardof}[1]{\mathopen{|\,}#1\mathclose{\,|}}
\newcommand{\ccardof}[1]{\mathopen{\,\|}#1\mathclose{\|\,}}
\newcommand{\ulsetof}[1]{\mathopen{\,|}#1\mathclose{|\,}}
\newcommand{\Card}{{\mathsf{Card}\/}}

\newcommand{\setof}[2]{\{#1\,:\,#2\}}
\newcommand{\ssetof}[1]{\{#1\}}

\newcommand{\subseteqand}[1]{\mathrel{\mathop{\subseteq}%
		\limits_{\scriptscriptstyle\hbox to 14pt{$\scriptscriptstyle #1$\hss}}}}

\newcommand{\dotcup}{\mathrel{\dot{\cup}}}

\newcommand{\mapping}[3]{#1:#2\rightarrow #3}
\newcommand{\isomrph}[3]{#1:#2\stackrel{\cong\hspace{0.8ex}}{\rightarrow}#3}
\newcommand{\isom}[3]{\isomrph{#1}{#2}{#3}}

\newcommand{\Elembed}[4]{#1:#2\stackrel{\prec\hspace{0.8ex}}{\rightarrow}_{#4}#3}
\newcommand{\EElembed}[4]{#1:#2\,\preccurlyeq_{#4}#3}
\newcommand{\eelembed}[3]{#1:#2\,\preccurlyeq\ #3}

\newcommand{\fnsp}[2]{\mbox{}^{{#1}\hspace{-0.02em}}#2}
\newcommand{\imageof}{{}^{\,{\prime}{\prime}}}
\newcommand{\numeralof}[1]{\ulcorner\hspace{-0.05ex}#1\hspace{-0.05ex}\urcorner}
\newcommand{\nnumeralof}[1]{%
  \ulcorner\ulcorner\hspace{-0.1ex}{#1}\hspace{-0.1ex}\urcorner\urcorner}

\newcommand{\seqof}[2]{\langle#1\,:\,#2\rangle}
\newcommand{\pairof}[1]{\langle#1\rangle}
\newcommand{\psof}[1]{{\mathcal P}\/(#1)}
\newcommand{\psnof}[2]{\ifx\bakakaba#1\bakakaba{\mathcal P}^n(#2)\else{\mathcal P}^#1(#2)\fi}

\newcommand{\modelof}[1]{\models\!\mbox{\rm``\,}#1\mbox{\,\rm''}}

\newcommand{\crit}{\mbox{\it crit\/}}

\newcommand{\circleq}{\mathrel{{\leqslant}%
		\hspace{-0.86ex}{\lower-0.53ex\hbox{$\scriptscriptstyle\circ$}}}}
\newcommand{\symb}[1]{{\mathord{\hspace{0.08em}\underbracket[0.6pt][2pt]{#1}}\hspace{0.08em}}}
\newcommand{\symbc}{\smash{\symb{c}}}
\newcommand{\symbin}{\mathbin{\symb{\in}}}

\newcommand{\restr}{\restriction}

\newcommand{\cf}{\mathop{cf\/}}
\newcommand{\Lim}{\mathop{Lim\/}}

\newcommand{\Col}{{\rm Col}}

\newcommand{\otp}{\mathop{\mbox{\it otp\/}}}

\newcommand{\dom}{\mathop{\it dom}}

\newcommand{\poP}{\bbd{P}}

\newcommand{\On}{{\rm On}}

\newcommand{\genG}{\mathbb{G}}

\newcommand{\LT}{{<}\,}
\newcommand{\LE}{{\leq}\,}
\newcommand{\GT}{{>}\,}

\newcommand{\llor}{{\bigvee\hspace{-1.5ex}\bigvee}\rule[-0.8ex]{0cm}{1ex}}
\newcommand{\lland}{{\bigwedge\hspace{-1.5ex}\bigwedge}\rule[-0.8ex]{0cm}{1ex}}
\ifarxived
\newcommand{\tenten}{\hspace{0.04ex}{.}\hspace{-0.04em}{.}\hspace{-0.04em}{.}\hspace{0.1ex}}
\newcommand{\ctenten}{,\mbox{}\hspace{0.08ex}{.}\hspace{-0.04em}{.}\hspace{-0.04em}{.}\hspace{0.1ex}}
\newcommand{\ctentenc}{,{}\linebreak[0]\hspace{0.04ex}{{.}\hspace{-0.04em}{.}\hspace{-0.04em}{.}\hspace{0.1ex},\,}\linebreak[0]}
\newcommand{\tentenc}{\hspace{0.04ex}{{.}\hspace{-0.04em}{.}\hspace{-0.04em}{.}\hspace{0.1ex},\,}\linebreak[0]}

\else
\newcommand{\tenten}{\hspace{0.04ex}{.}\hspace{-0.04zw}{.}\hspace{-0.04zw}{.}\hspace{0.1ex}}
\newcommand{\ctenten}{,\mbox{}\hspace{0.08ex}{.}\hspace{-0.04zw}{.}\hspace{-0.04zw}{.}\hspace{0.1ex}}
\newcommand{\ctentenc}{,{}\linebreak[0]\hspace{0.04ex}{{.}\hspace{-0.04zw}{.}\hspace{-0.04zw}{.}\hspace{0.1ex},\,}\linebreak[0]}
\newcommand{\tentenc}{\hspace{0.04ex}{{.}\hspace{-0.04zw}{.}\hspace{-0.04zw}{.}\hspace{0.1ex},\,}\linebreak[0]}

\fi
\newcommand{\xmbox}[1]{ $\relax{\rm #1}\relax$ }
\newcommand{\gmA}{\mathfrak{A}}
\newcommand{\gmB}{\mathfrak{B}}

\newcommand{\gmM}{\mathfrak{M}}

\newcommand{\continuum}{2^{\aleph_0}}

\newcommand{\calC}{{\mathcal C}}

\newcommand{\calF}{{\mathcal F}}

\newcommand{\calH}{{\mathcal H}}

\newcommand{\calL}{{\mathcal L}}

\newcommand{\calP}{{\mathcal P}}

\newcommand{\calS}{{\mathcal S}}

\newcommand{\calX}{{\mathcal X}}
\newcommand{\calY}{{\mathcal Y}}

\newcommand{\Lin}{{\calL}_{\in}}
\newcommand{\varin}{\mathrel{\varepsilon}}
\newcommand{\notvarin}{\mathrel{{\not}\,\varepsilon\,}}

\newcommand{\ZFC}{{\sf ZFC}}
\newcommand{\ZC}{{\sf ZC}}
\newcommand{\CH}{{\sf CH}}

\newcommand{\MA}{{\sf MA}}

\newcommand{\RP}{{\sf RP}}

\newcommand{\stat}{{stat}}
\newcommand{\Lstatinf}{\calL^{\aleph_0,\II}_{\stat,\kappa,\omega}}

\newcommand{\SDLS}{{\sf SDLS}}
\newcommand{\LSS}{{\sf LSTS}}
\newcommand{\II}{{\mathrm I\mathrm I}}
\newcommand{\III}{{\mathrm{III}}}
\newcommand{\HO}{{\mathrm H\mathrm O}}
\newcommand{\WCS}{{\sf WCS}}
\newcommand{\wcn}{{\mathfrak{wcn}}}
\newcommand{\cn}{{\mathfrak{cn}}}
\newcommand{\hn}{{\mathfrak{hn}}}
\newcommand{\lsn}{{\mathfrak{lstn}}}

\newcommand{\st}{such that}
\newcommand{\wrt}{with respect to}
\newcommand{\Wolog}{Without loss of generality}
\newcommand{\wolog}{without loss of generality}

\newcommand{\tfae}{the following are equivalent}

\newcommand{\po}{poset}
\newcommand{\pos}{posets}
\newcommand{\uniL}{{\sf L}}
\newcommand{\uniV}{{\sf V}}

\newcommand{\globalC}{{\sf C}}

\newcommand{\SINGL}{\mbox{$\mathsf{SING_L}$}}

\newcommand{\Pkl}[2]{\ifx\bakakaba#1\bakakababas\ifx\bakakaba#2\bakakaba{\mathcal 
    P}_\kappa(\lambda)\else{\mathcal P}_\kappa(#2)\fi\else{\mathcal P}_{#1}(#2)\fi}

\newcommand{\utildeT}[1]{%
  \hbox to 0pt{\smash{$\mathop{\textstyle #1}\limits_{%
			\raisebox{0.4ex}[0pt]{$\scriptstyle\sim$}}$}\hss}%
  \relax\phantom{\mathord{{#1}_{\rule[-0.6ex]{0pt}{1pt}}}}}
\newcommand{\utildeS}[1]{%
	\hbox to 0pt{\smash{$\mathop{\scriptstyle #1}\limits_{%
				\raisebox{0.6ex}[0pt]{$\scriptscriptstyle\sim$}}$}\hss}%
	\relax\phantom{\mathord{{#1}_{\rule[-0.6ex]{0pt}{1pt}}}}}
\newcommand{\utildeSS}[1]{%
	\hbox to 0pt{$\mathop{\scriptscriptstyle #1}%
		\limits_{\scriptscriptstyle\sim}$\hss}%
		\relax\phantom{\underline{#1}}}

{\end{minipage}\end{trivlist}}

\begin{document}
\maketitle
\renewcommand{\thefootnote}{$\ast$\ }
  \footnotetext{Graduate School of System Informatics, Kobe University \\Rokko-dai 1-1, Nada, Kobe 657-8501 Japan
   \\
    \quad\scalebox{0.95}[1]{\tt ${}^\dagger$ fuchino@diamond.kobe-u.ac.jp, ${}^\ddagger$ hsakai@people.kobe-u.ac.jp}}

\ifextended
\phantomsection
\addcontentsline{toc}{section}{* Weakly extendible cardinals and compactness of extended logics}
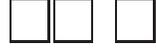
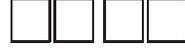
\addcontentsline{toc}{section}{**** by S.Fuchino and H.Sakai}
\fi

\addcontentsline{toc}{section}{Abstract}
\begin{abstract}
We introduce the notion of weakly extendible cardinals and show that these cardinals are characterized in 
terms of weak compactness of second order logic. The consistency strength and largeness of weakly 
extendible cardinals are located strictly between that of strongly unfoldable (i.e.\ shrewd) 
cardinals, and strongly uplifting cardinals. 

Weak compactness of many other logics can be connected to certain variants of the notion of 
weakly extendible cardinals.

We also show that, under $\uniV=\uniL$, a cardinal $\kappa$ is the weak compactness number 
of $\Lstatinf$ if and only if it is the weak compactness number 
of $\calL^\II_{\kappa,\omega}$. The latter condition is equivalent to the condition that $\kappa$ is weakly 
extendible by the characterization mentioned above (this equivalence holds without the assumption of $\uniV=\uniL$). 
\end{abstract}

\ifextended
{\extendedcolor 
\addcontentsline{toc}{section}{Contents}
\newcommand{\myscalebox}[1]{\scalebox{0.88}[1.06]{#1}}
\begin{quotation}
	\footnotesize
	\noindent
	\centerline{
      \normalsize\tt\quad\ Contents\hspace{6em}\mbox{}}\mbox{}\\
     {\mbox{}\hspace{-1.6em}\tt\makebox[3.4ex][l]{\ref{intro}.}%
      \hyperref[intro]{\tt\myscalebox{Compactness and L\"owenheim-Skolem theorems of 
          extended logics}}}\ \ \dotfill\ \ {\pageref{intro}}\\ 
      \noindent
    {\mbox{}\hspace{-1.6em}\tt\makebox[3.4ex][l]{\ref{weakl-ext}.}%
      \hyperref[weakl-ext]{\tt\myscalebox{Weakly extendible cardinals}}}\ \ \dotfill\ \ {\pageref{weakl-ext}}\\
      \noindent
    {\mbox{}\hspace{-1.6em}\tt\makebox[3.4ex][l]{\ref{comp-card}.}%
      \hyperref[comp-card]{\tt\myscalebox{General characterization of weakly $\calL$-compact cardinals}}}\ \ \dotfill\ \ {\pageref{comp-card}}\\
      \noindent
    {\mbox{}\hspace{-1.6em}\tt\makebox[3.4ex][l]{\ref{stat}.}%
      \hyperref[stat]{\tt\myscalebox{Weak compactness of stationary logic}}}\ \ \dotfill\ \ {\pageref{stat}}\\
      \noindent
     {\mbox{}\hspace{-1.6em}\hyperref[ref]{\tt
      References}}\ \ \dotfill\ \ {\pageref{ref}}\\ 

\end{quotation}}
\fi
\renewcommand{\thefootnote}{}
\footnotetext{{\it Date:} June 22, 2022
  \qquad {\it Last update:} 
  \today\ (\now\ JST)\vspace{-1\smallskipamount}
}
\footnotetext{{\it MSC2020 Mathematical Subject Classification:}
  03E35, 03E50, 03E55, 03E57, 03E65\vspace{-1\smallskipamount}}
\footnotetext{{\it Keywords:}
  L\"owenheim-Skolem theorems, compactness, extendible cardinals, subtle cardinals, strongly 
  uplifting cardinals, weakly extendible 
  cardinals, strongly unfoldable cardinals, shewed cardinals, strongly weakly 
  compact cardinals}
\footnotetext{\mbox{}\\[-2.52ex]The research is supported by Kakenhi Grant-in-Aid for Scientific Research (C) 20K03717}

\ifextended
\ifprivate
\footnotetext{This is a private extended version of a draft of the paper with the same title in 
  preparation.  
  All additional 
  details not to be contained the submitted version of the paper are either typeset in 
  {\extendedcolor "dark electric blue"} in case the text is also included in the extended 
  version of the paper or in {\darkgreen "dark green"} if the comment should appear only in 
  the private version. 
  The numbering of the assertions is kept identical with the submitted version.

  The most up-to-date file of this private extended version is downloadable as:\\
  {[0]} \,\href{https://fuchino.ddo.jp/papers/weakly-extendible-xx.pdf}{{\tt 
  https://fuchino.ddo.jp/papers/weakly-extendible-xx.pdf}} 
}
\else
\footnotetext{\extendedcolor This is an extended version of a draft of the paper with the same title.
  All additional 
  details not to be contained \extendedcolor in the submitted version of the paper are either typeset in 
  dark electric blue (the color in which this paragraph is typeset) or put in separate appendices. 
  The numbering of the assertions is kept identical with the submitted version.

  The most up-to-date file of the extended version is downloadable as:\\
  {[0]} \,\href{https://fuchino.ddo.jp/papers/weakly-extendible-x.pdf}{{\tt
      https://fuchino.ddo.jp/papers/weakly-extendible-x.pdf}} 
}
\fi
\else
\footnotetext{An updated and extended version of this paper with more details and 
  proofs is downloadable as:\\ {[0]}\ \,\href{https://fuchino.ddo.jp/weakly-extendible-x.pdf}{
   {\tt https://fuchino.ddo.jp/weakly-extendible-x.pdf}}}
\fi

\renewcommand{\thefootnote}{\arabic{footnote})\,}
\section{Compactness and L\"owenheim-Skolem theorems of extended logics}
\Label{intro} Suppose that $\calL$ is a logic with its model relation $\models_\calL$ (with all 
expected properties for a model relation, like those assumed in Lindstr\"om's 
theorem: see e.g.\ Chapter XII of \cite{ebbinghaus-etal}).  
In the following, we will not 
treat these logics and model relations associated to them in a strict axiomatic setting 
of an abstract model theory. Instead we just say that such logics 
are ``{\It proper}''\footnote{In \cite{MTL}, logics with corresponding ``natural'' properties are 
  called {\It regular}.} and, if necessary, 
mention only the  
specific properties of the logics explicitly which are assumed in some of the assertions.

For such {\It proper logic} $\calL$, a 
cardinal $\kappa$ is said to be {\It$\calL$-compact}, if for any collection $T$ 
of $\calL$-sentences, $T$ is  
{\It satisfiable} (i.e. there is a structure $\gmA$ with $\gmA\models_\calL T$) if (and only 
if) $T$ is {\It$\LT\kappa$-satisfiable} (i.e.\ all $T_0\in[T]^{\LT\kappa}$ are satisfiable).

If $\kappa$ is $\calL$-compact then any cardinal $\lambda>\kappa$ is $\calL$-compact as well. Thus 
by naming the minimal $\calL$-compact cardinal (if it exists), we completely describe the situation 
with $\calL$-compactness. We shall call this minimal cardinal the {\It compactness number} of $\calL$ and 
denote it by $\cn(\calL)$ (if it exists, otherwise we 
write $\cn(\calL)=\infty$). Thus
\begin{xitemize}
\item[] ${\darkred\cn(\calL)}:=\min(\setof{\kappa\in\Card}{\begin{array}[t]{@{}l}
    \mbox{for any }\calL\mbox{-theory }T\mbox{, }T\rlap{\mbox{ is satisfiable if and only}}\\
      \mbox{if all }T_0\in[T]^{\LT\kappa}\mbox{ are satisfiable}}\cup\ssetof{\infty}).
  \end{array}$
\end{xitemize}

We denote with  $\calL^\II$ the (monadic, full) second-order logic whose formulas are 
defined similarly to the first-order logic but with additional second-order variables $X$,
$Y$, $Z$ etc.\ and built-in predicate symbol $\varin$ for which it is allowed to build 
\begin{xitemize}
\xitem[x-intro-a-a] 
  atomic 
  formulas of the form $x\varin X$ for a first-order variable $x$ and a second-order 
  variable $X$: 
\end{xitemize}
the second-order variables are interpreted as such running over all subsets of the 
underlying set of the structure, and the atomic formula $x\varin X$ is interpreted as the (true) element relation between valuations 
of the variables (as element and subset of the underlying set of the structure in 
consideration). 

For a cardinal $\kappa$, {\darkred$\calL^\II_{\kappa,\omega}$} is the $\calL_{\kappa,\omega}$ 
extension of the second-order logic $\calL^\II$. That is, the formulas of
$\calL^\II_{\kappa,\omega}$ are constructed inductively as in $\calL^\II$ with the 
additional clause saying that conjunction $\lland\Phi$ and disjunction $\llor\Phi$ of 
set $\Phi$ of formulas of size $\LT\kappa$ are allowed as far as the set of free variables 
appearing in $\Phi$ is finite. The model relation for infinitary conjunction and 
disjunction is defined as expected. 

For $\omega<\lambda\leq\kappa$, {\darkred$\calL^\II_{\kappa,\lambda}$} is 
the $\calL_{\kappa,\lambda}$ extension of $\calL^\II_{\kappa,\omega}$ where we ease the 
restriction for $\lland\Phi$ and $\llor\Phi$ \st\ the formulas $\lland\Phi$ and $\llor\Phi$ are now permitted 
if the set of the free variables in $\Phi$ has cardinality $\LT\lambda$, and, in addition, 
existential and universal quantification over a block of quantifiers of size $\LT\lambda$ is 
allowed. 

For $n\geq 2$, let {\darkred$\calL^n$} be the $n$th-order logic and {\darkred$\calL^\HO$} be the higher-order 
logic defined as the union of all $\calL^n$, $n\in\omega$. {\darkred$\calL^\HO_{\kappa,\omega}$}, 
{\darkred$\calL^\HO_{\kappa,\lambda}$} are then defined similarly to
$\calL^\II_{\kappa,\omega}$, $\calL^\II_{\kappa,\lambda}$. 

Using this terminology, the classical characterization of extendible cardinals (for the 
definition of extendible cardinals see the beginning of the \sectionof{weakl-ext} below) by M.\ 
Magidor can be reformulated as follows:
\begin{Thm}{\rm (M.\,Magidor \cite{magidor}, see also Theorem 23.4 in 
    \cite{higher-inf})}\Label{P-intro-0}\wassertof{1}\ \  
  A cardinal $\kappa$ is extendible if and only if we have
  $\kappa=\cn(\calL^\II_{\kappa,\omega})$; if this equality 
  holds then we also have $\cn(\calL^\II_{\kappa,\omega})=\cn(\calL^\HO_{\kappa,\kappa})$. \smallskip

  \wassert{2} For a cardinal $\kappa$, $\kappa=\cn(\calL^\II)$ holds if and only if $\kappa$ is 
  the least extendible cardinal. If $\kappa$ is the least extendible cardinal then we also 
  have $\cn(\calL^\HO_{\kappa,\kappa})=\cn(\calL^\II_{\kappa,\omega})=\cn(\calL^\II_{\kappa,\kappa})
  =\cn(\calL^\II)=\kappa$. 
  \qed
\end{Thm}


\begin{Lemma}
  \Label{P-intro-4} For any uncountable cardinal $\kappa$, we have
  $\kappa\leq\cn(\calL_{\kappa,\omega})$. (Actually we have
  $\kappa\leq\wcn(\calL_{\kappa,\omega})$ where $\wcn(\cdot)$ is defined below). 
\end{Lemma}
\prf For all infinite $\mu<\kappa$, and infinite $\mu'\leq\mu$ the theory
\begin{xitemize}
\item[]   $T_{\mu'}:=\ssetof{\forall_x(\llor_{\alpha<\mu'}\,x\equiv c_\alpha)}\cup 
  \setof{d\not\equiv c_\alpha}{\alpha<\mu'}$. 
\end{xitemize}
is $\LT\mu'$-satisfiable but not satisfiable. This shows
$\mu<(\ \wcn(\calL_{\kappa,\omega})\,)\leq\cn(\calL_{\kappa,\omega})$.
\qedofLemma\qedskip

It is classical that strongly compact cardinals are also characterized in a similar vein:
\begin{Lemma}
  \Label{P-intro-1}
  An uncountable cardinal $\kappa$ is strongly compact if and only if
  $\kappa=\cn(\calL_{\kappa,\kappa})=\cn(\calL_{\kappa,\omega})$.
\end{Lemma}
\prf The first equality is simply the definition of strong compactness e.g.\ adopted in Kanamori 
\cite{higher-inf}. The second equation can be obtained by a modification of the proof of 
Proposition 4.1 in Kanamori \cite{higher-inf}. 
\qedofLemma
\qedskip

We define the 
{\It weak compactness 
spectrum} $\WCS(\calL)$ of a logic $\calL$ by: 
\begin{xitemize}
\item[] ${\darkred\WCS(\calL)}:=\setof{\kappa\in\Card}{{}
  \begin{array}[t]{@{}l}
    \kappa>\aleph_0,\,\mbox{for any }\calL\mbox{-theory }
    T\mbox{ of signature with }\\
    \mbox{at most }\LE\kappa\mbox{-many non-logical symbols, if }T\mbox{ is  }\\
    \LT\kappa\mbox{-satisfiable, then }T\mbox{ is satisfiable}
  }.  \end{array}$
\end{xitemize}

We shall also say that an uncountable cardinal $\kappa$ is {\It weakly $\calL$-compact}\/ if
$\kappa\in\WCS(\calL)$. 

In analogy to $\cn(\calL)$, we can also define the {\It weak compactness number}
$\wcn(\calL)$ as the minimum of $\WCS(\calL)$: 
\begin{xitemize}
\item[] \hspace{-0.4ex}${\darkred\wcn(\calL)}:=\min(\setof{\kappa\in\Card}{\begin{array}[t]{@{}l}
    \kappa>\aleph_0,\mbox{for any }\calL\mbox{-theory }T\mbox{ of signature with}\\ 
    \mbox{at most }\LE\!\kappa\mbox{-many non-logical symbols,\,if }T\mbox{ is }\\
    \LT\!\kappa\mbox{-satisfiable, then }T\mbox{ is satisfiable}}\cup\ssetof{\infty}).
  \end{array}$
\end{xitemize}
In contrast to $\calL$-compactness, there is no guarantee that $\WCS(\calL)$ is an 
end-segment of $\Card$ and hence $\wcn(\calL)$ does not necessarily decide $\WCS(\calL)$.

The following is a direct consequence of the definitions of the notions we introduced above: 

\begin{Lemma}
  \Label{P-intro-2}\wassertof{1}\ \ For any logic $\calL$, we have\
    $\setof{\kappa\in\Card}{\kappa\geq\cn(\calL)}\subseteq\WCS(\calL)$. \smallskip

  \wassert{2}  $\wcn(\calL)\leq\cn(\calL)$. \smallskip

  \wassert{3} If $\calL$ and $\calL'$ are logics \st\ each formulas of $\calL$ can be 
  translated to formulas of $\calL'$ then we have 
  \begin{xitemize}
  \item[] 
    $\wcn(\calL)\leq\wcn(\calL')$,
    $\cn(\calL)\leq\cn(\calL')$\ \  and\ \ \ $\WCS(\calL)\supseteq\WCS(\calL')$. 
    \qed
  \end{xitemize}
\end{Lemma}

\begin{Cor}
  \Label{P-intro-5} For cardinals $\kappa$ and $\kappa'$ with
  $\kappa\leq\kappa'$, if $\kappa'$ is the least extendible cardinal $\geq\kappa$, then
  \begin{xitemize}
  \item[] 
    $\setof{\lambda\in\Card}{\kappa'\leq\lambda}
    \subseteq\WCS(\calL^\HO_{\kappa,\kappa})\subseteq\WCS(\calL^\II_{\kappa,\omega})\subseteq\WCS(\calL^\II)$.
  \end{xitemize}

\end{Cor}
\prf By \Thmof{P-intro-0} and \Lemmaof{P-intro-2}.\qedofCor
\qedskip


\begin{Lemma}
  \Label{P-intro-3} \wassertof{1}\ \ An uncountable cardinal $\kappa$ is weakly compact if and only if
  $\kappa=\wcn(\calL_{\kappa,\kappa})=\wcn(\calL_{\kappa,\omega})$. 

  \wassert{2} $\setof{\mu\in\Card}{{}
    \begin{array}[t]{@{}l}
      \mu\geq\kappa, \mu\mbox{ is either weakly compact or }\\
      \mu\geq\mbox{the first strongly compact cardinal 
        above }\kappa}\ \subseteq\ \WCS(\calL_{\kappa,\kappa}).
    \end{array}
    $ 
\end{Lemma}
\memo{Does $=$ holds? Under $\uniV=\uniL$?}
\prf \assertof{1}: The first equality is just the definition of weakly compactness 
while the second equality can be obtained as a byproduct of the characterization of weak 
compactness by tree property (see e.g. the proof of Theorem 7.8 in Kanamori \cite{higher-inf} or 
Lemma 32.1 in Jech \cite{2nd-ed}). Note that by the definition of
$\wcn(\calL_{\kappa,\omega})$ (adopted from \cite{higher-inf}), it is easy to see that 
$\kappa=\wcn(\calL_{\kappa,\omega})$ implies that $\kappa$ is strongly inaccessible.
\ifprivate\memo{\privatecolor \mbox{Scan\char`_2022-04-20--14.47 L\"oskoT - annotated}, p.66}\fi
\smallskip

\assertof{2}: By \assertof{1}, \Lemmaof{P-intro-2}, and \Lemmaof{P-intro-1}. 
\qedofLemma\qedskip

The question about a possible solution of the following ``equation'' seems to be a natural 
one: 
\begin{xitemize}
\xitem[x-intro-a]
  $\displaystyle\frac{\mbox{ weakly compact cardinals }}{\mbox{ strongly compact cardinals}}
 \ \ =\ \ \frac{\mbox{$x$}}{\mbox{ extendible cardinals }}$.
\end{xitemize}

In terms of compactness and weak compactness, $x$ above must be a large cardinal property 
which should be characterized by 
\begin{xitemize}
\xitem[x-intro-a-0] 
  $\kappa=\wcn(\calL^\II_{\kappa,\omega})$. 
\end{xitemize}
In \sectionof{weakl-ext}, we introduce a new notion of large cardinals which we named ``weak  
extendibility'' and show that these cardinals are exactly those characterized by \xitemof{x-intro-a-0}
(\Thmof{P-w-ext-0}). 

We prove that the consistency strength of 
weak extendibility is strictly between that of subtleness and strong unfoldability 
(\Thmof{P-w-ext-1}). Note that, by L\"ucke \cite{luecke2}, strong unfoldability is 
equivalent to shrewdness of Michael Rathjen \cite{rathjen}.

Remembering Lindstr\"om's Theorem which gives the characterization of the first-order logic 
as the maximal logic satisfying the (countable) compactness and Downward L\"owenheim-Skolem 
Theorem down to countable (see e.g.\ \cite{ebbinghaus-etal}), it seems to be 
natural to consider also the spectrum of 
L\"owenheim-Skolem number of logics in our context (though a straight-forward 
generalization of Lindstr\"om's Theorem itself seems to be impossible: see \cite{shelah-vaananen}). 

In the following, we denote with $\darkred\ulsetof{\gmA}$ the underlying set of the (first-order) 
structure $\gmA$ and with $\darkred\ccardof{\gmA}$ the cardinality of (the underlying set of) $\gmA$.
Nevertheless, if we are talking about a set $A$, we continue to denote the cardinality of 
$A$ with $\cardof{A}$. 

Suppose that $\calL$ is a logic with the associated notion of elementary submodel $\prec_\calL$ (which 
should satisfy all the expected properties of an elementary submodel relation, like the 
properties assumed in the proof of the following \Lemmaof{P-intro-6}). The {\It 
  L\"owenheim-Skolem-Tarski spectrum} of $\calL$ is defined by: 
\begin{xitemize}
\item[] ${\darkred\LSS(\calL)}:=\setof{\mu\in\Card}{{}
  \begin{array}[t]{@{}l}
    \mbox{for any structure }\gmA\mbox{ of a countable signature}\\
    \mbox{and }S\subseteq\ulsetof{\gmA}\mbox{ with }\cardof{S}<\mu\mbox{, there is}\\
    \gmB\prec_\calL\gmA\mbox{ \st\ }S\subseteq\ulsetof{\gmB}
    \mbox{ and }\ccardof{\gmB}<\mu}.
  \end{array}$
  
\end{xitemize}

The terminology with ``L\"owenheim-Skolem-Tarski spectrum'' is chosen in accordance with 
that of \cite{magidor-vaananen}. The definition here is however slightly different since we 
consider not the cardinality of sentences but rather the cardinality of the signature. 
In \cite{sfetal-I}, and \cite{sfetal-II}, corresponding notion is refered to as  "strong L\"owenheim-Skolem property".

Our present definition of the L\"owenheim-Skolem-Tarski spectrum corresponds to the L\"owenheim-Skolem 
property in \cite{magidor}:  
\begin{Lemma}
  \Label{P-intro-6}
  For a logic $\calL$, we have
  \begin{xitemize}
  \item[] 
    $\LSS(\calL)=\setof{\mu\in\Card}{{}
        \begin{array}[t]{@{}l}
          \mbox{for any structure }\gmA\mbox{ with signature of}\\
          \mbox{size }\LT\mu\mbox{, there is }
          \gmB\prec_\calL\gmA\mbox{ \st\ }\ccardof{\gmB}<\mu}. 
        \end{array}$
  \end{xitemize}
\end{Lemma}
\prf
{\darkred ``$\darkred\subseteq$'':}\quad Suppose that $\darkred\mu\in\LSS(\calL)$ and let
$\gmA$ be a structure of 
signature of size $\nu<\mu$. \Wolog, we may assume that $\gmA$ is a relational structure 
and $\gmA=\pairof{\ulsetof{\gmA},R_{n,\alpha}}_{n\in\omega,\alpha<\nu}$ where 
$R_{n,\alpha}$ is an $n$-ary relation on $\ulsetof{\gmA}$ for $n\in\omega$ 
and $\alpha<\nu$. We may also assume, \wolog, that
$\ccardof{\gmA}\geq\mu$ and $\nu\subseteq\ulsetof{\gmA}$.\smallskip

Let ${\darkred R_n}:=\bigcup_{\alpha<\nu}\ssetof{\alpha}\times R_{n,\alpha}$ for each $n\in\omega$. 
Let ${\darkred \gmA^-}:=\pairof{\ulsetof{\gmA},R_n}_{n\in\omega}$.

Applying our assumption 
on $\mu$, we find $\gmB^-\prec_\calL\gmA^-$ with $\ccardof{\gmB^-}<\mu$ and
$\nu\subseteq\ulsetof{\gmB^-}$. By the last condition, we can reconstruct 
a submodel $\gmB$ of $\gmA$ from $\gmB^-$ with $\ulsetof{\gmB}=\ulsetof{\gmB^-}$
and $\gmB\prec_\calL\gmA$. 
\smallskip

{\darkred``$\darkred\supseteq$'':}\quad Suppose now that $\mu$ is in the set on the 
right side of the equality. Let 
$\gmA$ be a structure of size $\geq\mu$ with a countable  
signature, and $S\in[\,\ulsetof{\gmA}\,]^{\LT\mu}$.

Let $\gmA^+=\pairof{\gmA,a}_{a\in S}$. Applying the assumption on $\mu$, we obtain
$\gmB^+\prec_\calL\gmA^+$ of size $\LT\mu$. Denoting by $\gmB$ the structure $\gmB^+$ reduced to the 
original signature, we have $\ccardof{\gmB}<\mu$, $S\subseteq\ulsetof{\gmB}$ and
$\gmB\prec_\calL\gmA$. 
\qedofLemma

\begin{Lemma}
  \Label{P-intro-7}
  For any logic $\calL$, $\LSS(\calL)$ is a closed class of cardinals.
\end{Lemma}
\prf
Suppose that $\seqof{\kappa_\alpha}{\alpha<\delta}$ is a strictly increasing sequence 
in $\LSS(\calL)$ and $\kappa=\sup_{\alpha<\delta}\kappa_\alpha$. We want to show that
$\kappa\in\LSS(\calL)$.

Suppose that $\gmA$ is a structure of countable signature and $S\subseteq[\,\ulsetof{\gmA}\,]^{\LT\kappa}$. Let 
$\alpha<\delta$ be \st\ $\cardof{S}<\kappa_\alpha$. Since $\kappa_\alpha\in\LSS(\calL)$, 
there is a $\gmB\prec_\calL\gmA$ \st\ $S\subseteq\ulsetof{\gmB}$ and
$\ccardof{\gmB}<\kappa_\alpha<\kappa$. This shows that $\kappa\in\LSS(\calL)$.
\qedofLemma\qedskip

For $\calL^\II$, one detail of the definition of the $\calL^\II$-elementary submodel relation must be 
emphasized: for structures $\gmA$ and $\gmB$ with $\gmA\subseteq\gmB$, 
$\gmA\prec_{\calL^\II}\gmB$ if and only if for all $\calL^\II$-formula $\varphi(x_0\ctenten)$ (in the signature of 
the structures) \underline{without} \underline{free} \underline{second-order} 
\underline{variables} and $a_0\ctenten\in\ulsetof{\gmA}$,
\begin{xitemize}
\xitem[x-intro-0] $\gmA\models_{\calL^\II}\varphi(a_0\ctenten)$\ \ $\Leftrightarrow$\ \ 
  $\gmB\models_{\calL^\II}\varphi(a_0\ctenten)$.
\end{xitemize}

Exclusion of the second-order parameters from the definition is justified by the fact that 
we would have only trivial cases of elementary substructure relation
$\gmA\prec_{\calL^\II}\gmB$ (namely when $\gmA=\gmB$) if we would have included the 
second-order parameters in the definition \xitemof{x-intro-0}. 

By \Lemmaabove, the proof of Theorem 1 and Theorem 2 in \cite{magidor} can be recast to 
show the following:
\begin{Thm}{\rm (M.\,Magidor \cite{magidor})}\Label{P-intro-7-0}\smallskip

  $\LSS(\calL^\II)=\LSS(\calL^\HO)=\setof{\kappa\in\Card}{{}
    \begin{array}[t]{@{}l}
      \kappa\mbox{ is a supercompact cardinal,}\\
      \mbox{or a 
        limit of supercompact cardinals}}.
    \end{array}$\vspace{-1.2em}\\
\qed
\end{Thm}
\imemo{In extended version, add the proof from kobe-fuchino-2022-06-01.tex}

Similarly to $\cn(\calL)$ and $\wcn(\calL)$, we define the L\"owenheim-Skolem-Tarski number of a 
logic $\calL$ be the least element of $\LSS(\calL)$. More precisely, for a logic $\calL$, we let
\begin{xitemize}
\xitem[x-intro-1-0] 
  $\lsn(\calL):=\min(\LSS(\calL)\cup\ssetof{\infty})$.
\end{xitemize}

\Thmof{P-intro-7-0} implies 
\begin{xitemize}
\xitem[x-intro-1-1] 
  $\lsn(\calL^\II)=\lsn(\calL^\HO)=$ the least 
  supercompact cardinal (if there is one). 
\end{xitemize}

$\calL^{\aleph_0,\II}$ denotes the weak (monadic) second-order logic with second-order 
variables $X$, $Y$, $Z$ etc. whose intended interpretation is that they run over countable 
subsets of the underlying set of the structure in consideration. We shall call this type of 
second-order variables weak second-order variables (in $\aleph_0$ -interpretation).

The formulas of $\calL^{\aleph_0,\II}$, $\calL^{\aleph_0,\II}_{\kappa,\omega}$,
$\calL^{\aleph_0,\II}_{\kappa,\lambda}$ are defined in exactly the same way as the formulas 
of  $\calL^{\II}$, $\calL^{\II}_{\kappa,\omega}$, 
$\calL^{\II}_{\kappa,\lambda}$ but the 
inductive definition of the semantics uses the $\aleph_0$-interpretation: 

\begin{xitemize}
\item[] $\gmA\models_{\calL^{\aleph_0,\II}}\exists X\varphi(a_0\ctentenc X,A_0\ctenten)$\\
  \qquad$:\Leftrightarrow$\ \ 
  there is $B\in[\,\ulsetof{\gmA}\,]^{\aleph_0}$ \st\
  $\gmA\models_{\calL^{\aleph_0,\II}}\varphi(a_0\ctentenc B,A_0\ctenten)$. 
\end{xitemize}

In case of this weak second-order logic, a definition of the elementary submodel relation with second-order 
parameters also makes sense.
Here, for simplicity, we assume below that the elementarity
$\prec_{\calL^{\aleph_0,\II}}$, $\prec_{\calL^{\aleph_0,\II}_{\kappa,\omega}}$ etc.\ are 
always defined similarly to \xitemof{x-intro-0} without second-order parameters. 

The L\"owenheim-Skolem theorems of stationary quantifier in the context of the weak 
second-order logic and its interplay with various reflection principles were studied in 
\cite{sfetal-I}. 

The logic $\calL^{\aleph_0}_\stat $ is defined as the logic with monadic second-order 
variables with the second-order quantifier $\stat$ (and its dual $a.a.$) \st\ the 
recursive definition of $\calL^{\aleph_0}_\stat $-formulas includes \xitemof{x-intro-a-a} 
together with the following in 
addition to the clauses in the usual definition of the first-order logic: 
\begin{xitemize}
\item[] If $\varphi$ is a $\calL^{\aleph_0}_\stat $-formula and $X$ a second-order variable, then 
  $\stat X\,\varphi$ is also a $\calL^{\aleph_0}_\stat $-formula.
\end{xitemize}
The model relation of this logic is then defined as usual with the following additional clause:
\begin{xitemize}
\item[] For any $\calL^{\aleph_0}_\stat $-formula
  $\varphi(x_0\ctentenc X, X_0\ctenten)$, 
  $a_0\ctenten\in\ulsetof{\gmA}$ and
  $A_0\ctenten\in[\,\ulsetof{\gmA}\,]^{\aleph_0}$,
  \begin{xitemize}
  \item[] 
    $\gmA\models_{\calL^{\aleph_0}_\stat }stat\, X\ \varphi(a_0\ctentenc X, A_0\ctenten)$ 
    \\[2\jot]
    :$\Leftrightarrow$\ \ $\setof{A\in[\,\ulsetof{\gmA}\,]^{\aleph_0}}{
      \gmA\models_{\calL^{\aleph_0}_\stat}\varphi(a_0\ctentenc A, A_0\ctenten)}$ 
    is a stationary\\
    \phantom{:$\Leftrightarrow$\ \ } subset of $[\,\ulsetof{\gmA}\,]^{\aleph_0}$. 
  \end{xitemize}
\end{xitemize}

$\calL^{\aleph_0,\II}_\stat$ is the extension of $\calL^{\aleph_0}_\stat$ with 
the weak second-order existential (and universal) quantification. 
$\calL^{\aleph_0}_{stat, \kappa,\lambda}$ and $\calL^{\aleph_0,\II}_{stat, \kappa,\lambda}$ 
are then the infinitary versions of $\calL^{\aleph_0}_\stat$ and
$\calL^{\aleph_0,\II}_\stat$ defined as expected.

\begin{Lemma}
  \Label{P-intro-8} \wassertof{1}\ \  The expressive power of
  $\calL^{\aleph_0,\II}_{\kappa,\omega}$ exceeds that of $\calL_{\kappa,\omega}$. 
  \smallskip\memo{\darkred and strictly stronger than 
    $\calL^{\aleph_0,\II}$?}

  \wassert{2} For $\omega<\lambda\leq\kappa$, $\calL^{\aleph_0,\II}_{\kappa,\lambda}$ has 
  the same expressive power as $\calL_{\kappa,\lambda}$. \smallskip

  \wassert{3} For 
  any $\omega<\lambda\leq\kappa$, $\calL^{\aleph_0,\II}_{stat,\kappa,\lambda}$ is 
  interpretable in $\calL^\III_{\kappa,\lambda}$ where $\calL^\III$ denotes the 
  third-order logic, and $\calL^\III_{\kappa,\lambda}$ its infinitary extension. 
\end{Lemma}
\prf \assertof{1}: For a binary relation symbol $R$,  ``$R$ is well-founded'' can be expressed by the
$\calL^{\aleph_0,\II}$-sentence:
\begin{xitemize}
\xitem[] $\forall X\exists_y(y\varin X\ \land\ \forall_x(x\not\mathrel{R}y))$.
\end{xitemize}
On the other hand, $\calL_{\kappa,\omega}$ cannot express the well-foundedness of $R$ by a 
theorem of Lopez-Escobar \cite{lopez-escobar}. \smallskip

\assertof{2}: We define the translation of $\calL^{\aleph_0,\II}_{\kappa,\lambda}$-formula $\varphi$ into
$\calL_{\kappa,\lambda}$-formula $\varphi_*$ by assigning each second-order variable $X$ to 
countably many new first order variables $x^X_i$, $i\in\omega$; assigning each atomic 
formula of the form $x\varin X$ to the formula 
\begin{xitemize}
\item[] 
  $\llor\setof{x\equiv x^X_i}{i\in\omega}$; 
\end{xitemize}
and assigning each $\calL^{\aleph_0,\II}_{\kappa,\lambda}$-formula $\varphi$ of the form
$\exists X\,\psi(\tentenc,X\ctenten)$  
to $\calL_{\kappa,\lambda}$-formula $\varphi_*$ of the form 
$\exists_{\,x^X_0}\exists_{\,x^X_1}\tenten(\,\lland_{i<j<\omega}x^X_i\not\equiv x^X_j\ \land\ 
\psi_*(\tentenc x^X_0,x^X_1\ctentenc\tenten))$. \smallskip

\assertof{3}: A translation of $\calL^{\aleph_0,\II}_{stat,\kappa,\lambda}$-formula $\varphi$ 
into $\calL^\III_{\kappa,\lambda}$-formula $\varphi_{\ast\ast}$ will do which starts 
similarly to $\varphi\mapsto\varphi_*$ in \assertof{2} and continues with 
the following details of the recursive definition: 

If the $\calL^{\aleph_0,\II}_{stat,\kappa,\lambda}$-formula $\varphi$ is of the form
$\exists X\psi$, the translation $\varphi_{\ast\ast}$ is defined to be the $\calL^\III_{\kappa,\lambda}$-formula:
\begin{xitemize}
\item[] 
 $\exists X\exists_{\,x^X_0}\exists_{\,x^X_1}\cdots(\lland_{i<j<\omega}x^X_i\not\equiv x^X_j\ \land\ 
  \forall_y(y\varin X\ \leftrightarrow\ \llor_{k<\omega}y\equiv x^X_k)\ \land\ \psi_{\ast\ast})$.
\end{xitemize}

If the $\calL^{\aleph_0,\II}_{stat,\kappa,\lambda}$-formula $\varphi$ is of the form 
$stat\, X\psi$, the translation $\varphi_{\ast\ast}$ is defined to be the $\calL^\III_{\kappa,\lambda}$-formula:
\begin{xitemize}
\item[] 
  $\exists\calX\,(\ \forall X\,(X\varin\calX\ \rightarrow\\
  \qquad\quad\ \ \ 
  \exists_{\,x^X_0}\exists_{\,x^X_1}\tenten(\lland_{i<j<\omega}x^X_i\not\equiv x^X_j\ \land\ 
  \forall_y(y\varin X\ \leftrightarrow\ \llor_{k<\omega}y\equiv x^X_k)\ \land\ \psi_{\ast\ast}))\\[\jot]
  \quad\quad
  \land\ \forall\calY\,(\mbox{``}\calY\mbox{ is a club of countable sets''}
  \ \rightarrow \exists X\,(X\varin\calX\ \land\ X\varin\calY))
)$
\end{xitemize}
where $\calX$ and $\calY$ are third order variables. 
\qedofLemma

\begin{Prop}
  \Label{P-intro-9} \wassertof{1}\ \ An uncountable cardinal $\kappa$ is weakly compact if and only if
  $\kappa=\wcn(\calL^{\aleph_0,\II}_{\kappa,\omega})$.\smallskip

  \wassert{2}
  $\wcn(\calL^{\aleph_0,\II}_{stat,\kappa,\omega})\leq\wcn(\calL^\III_{\kappa,\omega})
  \leq\mbox{the first weakly extendible cardinal above }\kappa$.
\end{Prop}
\prf By \Lemmaof{P-intro-3}, \Lemmaof{P-intro-2},\,\assertof{3} and \Lemmaof{P-intro-8}.
The rightmost inequality in \assertof{2} follows from \Thmof{P-w-ext-0} (for the definition 
of weakly extendibility, see around \xitemof{x-w-ext-0}). 
\qedofProp
\qedskip

In \sectionof{stat}, we discuss about further results related to \Propof{P-intro-9} 
above. In particular, we shall show that $\wcn(\calL^{\aleph_0,\II})$ is above certain 
large cardinals (see the remark after \Propof{P-stat-0}) while, in $L$, the condition
$\kappa=\wcn((\calL^{\aleph_0}_\stat)_{\kappa,\omega})$ is equivalent to
$\kappa=\wcn((\calL^\II)_{\kappa,\omega})$ (\Thmof{P-stat-7}).

\section{Weakly extendible cardinals}
\Label{weakl-ext}
Let us begin with recalling the definition of extendible cardinals: 
a cardinal $\kappa$ is {\It extendible} if, for any $\eta>0$, there is a $\zeta$ and $j$ 
\st\ $\Elembed{j}{V_{\kappa+\eta}}{V_\zeta}{\kappa}$ (see Kanamori 
\cite{higher-inf})\footnote{\Label{fn-0}In the following, write 
  $\Elembed{j}{M}{N}{\kappa}$ to denote the situation that $j$ is an elementary 
  embedding (in the sense of elementarity in the first-order logic) of $M=\pairof{M,\in}$ 
  into $N=\pairof{N,\in}$, both $M$  and $N$  are transitive, and $\kappa$ is the critical point 
  of $j$. $M$, $N$, $j$ are sets in the definition of extendibility but later we shall also treat the cases 
  where these are proper classes possibly in some generic extension of $\uniV$.}.

As it is well known, weakly and strongly compact cardinals are characterized in terms of elementary embeddings:
\begin{Lemma}
  \Label{P-w-ext-a}
 \wassertof{1}\ {\rm (see e.g.\ \cite{joel})}\ A 
  cardinal $\kappa$ is weakly compact if and only  
  if $\kappa=2^{\LT\kappa}$ and satisfies the following {\It Embedding Property}: for any 
  transitive $M$  
  with $\kappa\in M$ and $\cardof{M}=\kappa$, there is a 
  transitive $N$ with j \st\ $\Elembed{j}{M}{N}{\kappa}$. \smallskip

  \wassert{2} {\rm(Theorem 22.17 in \cite{higher-inf})}\ A cardinal $\kappa$ is 
  strongly compact if and only if, for any cardinal 
  $\lambda>\kappa$, there are classes $M$, $j$ \st\ $\Elembed{j}{\uniV}{M}{\kappa}$ and, 
  for any $a\in[M]^{\leq\lambda}$, there is $b\in[M]^{\LT j(\kappa)}\cap M$ with
  $a\subseteq b$. \qed 
\end{Lemma}

Comparing the definition of extendibility with the characterizations of weak and strong 
compact cardinals in \Lemmaabove, The following notion seems to be a good candidate of the 
large cardinal property which should be characterized by \xitemof{x-intro-a-0}:

We shall say that a cardinal $\kappa$ is {\It weakly extendible} if 
\begin{xitemize}
\xitem[x-w-ext-0] 
  $\kappa=2^{\LT\kappa}$ 
  and,  
\xitem[x-w-ext-1] 
  for any
  $\theta>\kappa$, and $M\prec V_\theta$ with $\kappa+1\subseteq M$, and $\cardof{M}=\kappa$, 
  there are $\overline{\theta}$ and $j$ with $\EElembed{j}{M}{V_{\overline{\theta}}}{\kappa}$.
\end{xitemize}

Here, when we write $\EElembed{j}{M}{N}{\kappa}$, we assume that $\kappa+1\subseteq M, N$ (but 
$M=\pairof{M,\in}$ and $N=\pairof{N,\in}$ are not necessarily transitive), and $\mapping{j}{M}{N}$ is an 
elementary embedding with $j\restr\kappa=id_\kappa$ and $j(\kappa)>\kappa$. Similarly to 
the notation  $\Elembed{j}{M}{N}{\kappa}$, the ``elementary 
embedding'' here is meant in terms of first-order logic.  
Note that, if
$\isom{m}{M}{M_0}$, $\isom{n}{N}{N_0}$ are the Mostowski collapses of $M$ and $N$, we have
$\Elembed{n\circ j\circ m^{-1}}{M_0}{N_0}{\kappa}$. 

A weakly extendible cardinal is (strongly) inaccessible. This can be shown by an easy direct 
argument\ifextended{ \extendedcolor (see \LemmaAof{P-w-ext-a-0} below)}\fi. However, by the definition of the weak 
extendibility and the characterization of weak compactness 
(\Lemmaof{P-w-ext-a},\,\assertof{1}), 
we see immediately that a weakly extendible cardinal is weakly compact and we know that a weakly compact 
cardinal is inaccessible (e.g.\ Proposition 4.4 in \cite{higher-inf}).

\ifextended{\extendedcolor
\begin{LemmaA}
  \Label{P-w-ext-a-0} If $\kappa$ is weakly extendible then $\kappa$ is (strongly) inaccessible.
\end{LemmaA}
\prf $\kappa$ is not a successor cardinal: Suppose otherwise, say $\kappa=\mu^+$. Let
$M\prec V_{\kappa+1}$ be \st\ $\kappa+1\subseteq M$ and $\cardof{M}=\kappa$. Let
$\EElembed{j}{M}{V_{\overline{\theta}}}{\kappa}$. Then $j(\mu)=\mu$ and hence
$V_{\overline{\theta}}\modelof{j(\kappa)\mbox{ is the successor of }\mu}$. It follows that
$j(\kappa)=(\mu^+)^{V_{\overline{\theta}}}\leq\mu^+=\kappa$. This is a contradiction.
\smallskip

$\kappa$ is not a singular cardinal: Suppose otherwise, say 
$\kappa=\lim_{\alpha<\mu}\kappa_\alpha$ for $\mu<\kappa$ and $\kappa_\alpha<\kappa$ for
$\alpha<\mu$. For an $M\prec V_{\kappa+1}$ with $\seqof{\kappa_\alpha}{\alpha<\mu}\in M$ 
and $\cardof{M}=\kappa$ let $\EElembed{j}{M}{V_{\overline{\theta}}}{\kappa}$. By 
elementarity, we have
\begin{xitemize}
\item[] $V_{\overline{\theta}}\modelof{j(\kappa)\mbox{ is the limit of }j(\seqof{\kappa_\alpha}{\alpha<\mu})}$.
\end{xitemize}
Since $j(\seqof{\kappa_\alpha}{\alpha<\mu})=\seqof{\kappa_\alpha}{\alpha<\mu}$, it follows 
that $j(\kappa)=\kappa$. This is again a contradiction.  
\smallskip

$\kappa$ is a strong limit: Otherwise there is a $\mu<\kappa$ \st\ $2^\mu=\kappa$ (Note 
that the assumption $2^{\LT\kappa}=\kappa$ implies $2^\mu\leq\kappa$).

Let $\mapping{f}{\fnsp{\mu}{2}}{\kappa}$ be a bijection and let $M\prec V_{\kappa+1}$ be 
\st\ $f\in M$, $\kappa+1\subseteq M$ and $\cardof{M}=\kappa$. 
Note that $\fnsp{\mu}{2}\subseteq M$ and hence $(\fnsp{\mu}{2})^M=\fnsp{\mu}{2}$. 

Let
$\EElembed{j}{M}{V_{\overline{\theta}}}{\kappa}$. Then $j(\fnsp{\mu}{2})=\fnsp{\mu}{2}$ and 
$j(f)=j(f)\restr\fnsp{\mu}{2}=f$. By elementarity, it follows that
$V_{\overline{\theta}}\modelof{j(\kappa)\mbox{ is the least upper bound of }f\imageof\,{\fnsp{\mu}{2}}}$ 
which implies $j(\kappa)=\kappa$. This is a contradiction. \qedofLemma\qedskip}\fi 

The next theorem shows that our notion of weak extendibility is exactly what we are looking for. 

After we had written the first version of this paper, Will Boney told us that he also 
obtained a slight variation of the following theorem which also can be seen as a characterization of the weak 
compactness of infinitary second-order logics (cf.\ Theorem 4.5 in \cite{boney}). We keep 
our proof of the theorem here since it will be modified to obtain further results in this 
and next sections. 

\begin{Thm}
  \Label{P-w-ext-0} For a cardinal $\kappa$, \tfae:\smallskip

  \wassert{a} $\kappa=\wcn(\calL^\II_{\kappa,\omega})$. 
  
  \wassert{b} $\kappa=\wcn(\calL^\HO_{\kappa,\kappa})$.

  \wassert{c} $\kappa$ is weakly extendible. 

\end{Thm}
\prf ``\assertof{c} $\Rightarrow$ \assertof{b}'': Assume that $\kappa$ is weakly 
extendible. $\kappa$ is then inaccessible (see the remark after the definition 
of 
weakly extendibility). Suppose that $T$ is a
$\LT\kappa$-satisfiable $\calL^\HO_{\kappa,\kappa}$-theory of signature of size
$\leq\kappa$. We want to show that $T$ is satisfiable.

Since $T$ has cardinality $\leq\kappa$, we may assume that $T$ is a subset of $\kappa$ by some 
reasonable coding. 

Let $\theta$ be large enough. In particular, \st
\begin{xitemize}
\xitem[x-w-ext-2] $\theta\geq\kappa^{+\omega}$ and 
\xitem[x-w-ext-3] 
  $V_\theta\modelof{\,T\mbox{ is }\LT\kappa\mbox{-satisfiable}}$.
\end{xitemize}
Let $M\prec V_{\theta}$ be \st\ $\kappa+1\subseteq M$, $T\in M$, $\cardof{M}=\kappa$,  
and let $\overline{\theta}$, $j$ be \st\
$\EElembed{j}{M}{V_{\overline{\theta}}}{\kappa}$. Then we have
$V_{\overline{\theta}}\modelof{j(T)\mbox{ is }\LT j(\kappa)\mbox{-satisfiable}}$ by 
\xitemof{x-w-ext-3} and by elementarity of $j$.

Since 
$V_{\overline{\theta}}\modelof{\cardof{T}<j(\kappa)\mbox{ and }
\smash{\obecause{T}{0pt}{\qquad\qquad\quad$=j(T)\cap\kappa$}}\!\!\subseteq j(T)}$, it follows that there is
$\gmA\in V_{\overline{\theta}}$ \st\
$V_{\overline{\theta}}\modelof{\gmA\models_{\calL^\HO_{\kappa,\kappa}}T}$. 
Now $\overline{\theta}\geq j(\kappa)^{+\omega}\geq \kappa^{+\omega}$ by \xitemof{x-w-ext-2} and by elementarity 
of $j$. Thus, it follows that $\gmA\models_{\calL^\HO_{\kappa,\kappa}}T$. Thus, 
$T$ is realizable.\smallskip

``\assertof{b} $\Rightarrow$ \assertof{a}'': By \Lemmaof{P-intro-2},\,\assertof{3} and  
\Lemmaof{P-intro-4}. \smallskip

``\assertof{a} $\Rightarrow$ \assertof{c}'': Assume that
$\kappa=\wcn(\calL^\II_{\kappa,\omega})$ holds. Then we 
have $\kappa=\wcn(\calL_{\kappa,\omega})$ by \Lemmaof{P-intro-4}. Hence $\kappa$ is weakly 
compact by \Lemmaof{P-intro-3}. In particular, $\kappa$ is inaccessible. Thus it is enough 
to show that $\kappa$ satisfies \xitemof{x-w-ext-1}.

Suppose that $\theta>\kappa$ and $M\prec V_\theta$ is \st\ $\kappa+1\subseteq M$ and
$\cardof{M}=\kappa$.

Let $\darkred\varphi^*$ be an $\calL^\II$-sentence in the signature $\ssetof{\symb{\in}}$ \st\ 
\begin{xitemize}
\xitem[x-w-ext-3-0] $\pairof{\ulsetof{\gmA},\symb{\in}^\gmA}\models\varphi^*$\ \ $\Leftrightarrow$\ \ 
  $\symb{\in}^\gmA$ is well-founded and 
  extensional binary relation,\\[\jot]
  \phantom{$\pairof{\ulsetof{\gmA},\symb{\in}^\gmA}\models\varphi^*$\ \ $\Leftrightarrow$\ \ }%
  and the Mostowski collapse of $(\pairof{\ulsetof{\gmA},\symb{\in}^\gmA})$ is $\pairof{V_\gamma,\in}$\\ 
  \phantom{$\pairof{\ulsetof{\gmA},\symb{\in}^\gmA}\models\varphi^*$\ \ $\Leftrightarrow$\ \ }%
  for some $\gamma$
\end{xitemize}\smallskip

Let
\begin{xitemize}
\xitem[x-w-ext-4]
  ${\darkred T}:=\ssetof{\varphi^*}\ 
  \begin{array}[t]{@{}l}
    \cup\ \setof{\varphi(\symbc_{a_0}\ctenten)}{{}
      \begin{array}[t]{@{}l}
        \varphi
        \mbox{ is a first-order formula in the signature }\ssetof{\symb{\in}},\\
        a_0\ctenten\in M\mbox{ and }M\models\varphi(a_0\ctenten)}
      \end{array}
      \\[\jot]
    \cup\ \setof{\forall_x(x\mathbin{\symb{\in}}\symbc_\alpha\ \leftrightarrow\ 
      \llor_{\beta<\alpha}x\equiv\symbc_\beta)}{\alpha<\kappa}\\
    \cup\ \setof{\symbc_\alpha\mathbin{\symb{\in}}\symb{d}}{\alpha<\kappa}\\
    \cup\ \ssetof{\symb{d}\mathbin{\symb{\in}}\symbc_\kappa}.\qquad\qquad\qquad\qquad
    \qquad\qquad\qquad\qquad\qquad
    \ifarxived 
    \smash{\raisebox{2.88ex}{$\left.\rule[0pt]{0pt}{4.68ex}\right\}$\quad  \assertof{*}}}
    \else
    \smash{\raisebox{1.6zh}{$\left.\rule[0pt]{0pt}{2.6zh}\right\}$\quad  \assertof{*}}}
    \fi
  \end{array}
  $
\end{xitemize}
The signature of the $\calL^\II_{\kappa,\omega}$-theory $T$ is
$\ssetof{\symb{\in}, \symb{d}}\cup\setof{\symbc_a}{a\in M}$ and it is 
of cardinality $\kappa$. 
\begin{Claim}
  $T$ is $\LT\kappa$-satisfiable. 
\end{Claim}
\prfofClaim Suppose that $T_0\in[T]^{\LT\kappa}$. We have to show that $T_0$ has a model. 

Let $C_0:=\setof{a\in M}{\symbc_a\mbox{ appears in }T_0}$ and 
\begin{xitemize}
\item[] 
  $\alpha^*:=\sup\setof{\alpha+1}{\alpha<\kappa,\, \symbc_\alpha\xmbox{ appears in }T_0}$. 
\end{xitemize}
Let $\gmM=\pairof{V_\theta,\symb{d}^\gmM,(\symbc_a)^\gmM,\symb{\in}^\gmM}_{a\in C_0}$ where 
$\symb{d}^\gmM:=\alpha^*$ and $(\symbc_a)^\gmM:=a$ for all $a\in C_0$.
Then $\gmM\models_{\calL^\II_{\kappa,\omega}}T_0$. 
\qedofClaim\qedskip

By 
the assumption on $\kappa$, it follows that $T$ is satisfiable.
Let $\gmB$ be a model of $T$. By $\gmB\models_{\calL^\II}\varphi^*$, we can take the Mostowski collapse 
$\gmB^*$ of $\gmB$, and 
$\ulsetof{\gmB^*}=V_{\overline{\theta}}$ for some ordinal $\overline{\theta}$. 
Note that we have $\symb{\in}^{\gmB^*}=\,\in$. By the definition $\mapping{j}{M}{V_{\overline{\theta}}}$\,; 
$a\mapsto [\symbc_a]^{\gmB^*}$, we obtain $\EElembed{j}{M}{V_{\overline{\theta}}}{\kappa}$. 
Note that $\kappa$ is the critical point of $j$ by \xitemof{x-w-ext-4}\,\assertof{*}. 
\qedofThm\qedskip

We also obtain a characterization of elements of $\WCS(\calL^\II)$ by modifying the proof 
of \Thmof{P-w-ext-0}.

\memo{Later in a more general setting the 
  condition ``$j(\kappa)>\sup(j\imageof\kappa)$'' is really needed since $j\imageof{T}$ 
  might not be in the target model.} 
Let us call a cardinal $\kappa$ {\It weakly sub-extendible} if, for any $\theta>\kappa$ and
$M\prec V_\theta$ with $\kappa+1\subseteq M$ and $\cardof{M}=\kappa$, there are 
$\overline{\theta}$ and $j$ with $\EElembed{j}{M}{V_{\overline\theta}}{\mu}$ for some $\mu\leq\kappa$ 
\st\ \ixitem[x-w-ext-4-0] $j(\kappa)>\sup(j\imageof\kappa)$. 


\begin{Thm}
  \Label{P-w-ext-0-0}
  For a cardinal $\kappa$, $\kappa\in\WCS(\calL^\II)$\ \ $\Leftrightarrow$\ \ $\kappa$ is 
  weakly sub-extendible.
\end{Thm}
\prf $(\Leftarrow)$:\quad Suppose that $\kappa$ is weakly sub-extendible and $T$ is 
a $\LT\kappa$-satisfiable $\calL^\II$-theory of cardinality $\kappa$ (note that, since $\calL^\II$ is finitary, 
the cardinality of $T$ ($+\aleph_0$) is equal to the cardinality of the signature
($+\aleph_0$). We may assume that $T$ is is nicely coded as a subset of $\kappa$. 

Let $\theta>\kappa$ be large enough regular cardinal \st\
$V_\theta\modelof{T\mbox{ is }\LT\kappa\mbox{-satisfiable}}$. Let $M\prec V_\theta$ be \st\
$\kappa+1\subseteq M$, $T\in M$ and $\cardof{M}=\kappa$. By assumption there are $j$ and
$\overline{\theta}$ \st\ $\EElembed{j}{M}{V_{\overline{\theta}}}{\mu}$ for some 
$\mu\leq\kappa$ and $j(\kappa)>\kappa$.

By elementarity, we have
$V_{\overline{\theta}}\modelof{j(T)\mbox{ is }\LT j(\kappa)\mbox{-satisfiable}}$. 
Since $V_{\overline{\theta}}\modelof{j\imageof T\in[j(T)]^{\LT j(\kappa)}}$, it follows 
that $V_{\overline{\theta}}\modelof{j\imageof T\mbox{ is satisfiable}}$. Let 
$\gmA\in V_{\overline{\theta}}$ be \st\
$V_{\overline{\theta}}\modelof{j\imageof T\models_{\calL^\II}\gmA}$. Since
$V_{\overline{\theta}}$ interpret $\models_{\calL^\II}$ correctly (later we introduce the 
terminology with which we would say $V_\theta$ is $\calL^\II$-truthful), we have
$\gmA\models_{\calL^\II}j\imageof{T}$. Thus, by renaming the structure of $\gmA$ we obtain 
a model $\gmA'$ of $T$. \smallskip

$(\Rightarrow)$:\quad Suppose that $\kappa\in\WCS(\calL^\II)$. 
Suppose that $\theta$ and $M$ are as in the definition of the weakly sub-extendibility. Let
\begin{xitemize}
\xitem[x-w-ext-5]
  ${\darkred T}:=\ssetof{\varphi^*}\ 
  \begin{array}[t]{@{}l}
    \cup\ \setof{\varphi(\symbc_{a_0}\ctenten)}{{}
      \begin{array}[t]{@{}l}
        \varphi
        \mbox{ is a first-order formula in the signature }\ssetof{\symb{\in}},\\
        a_0\ctenten\in M\mbox{ and }M\models\varphi(a_0\ctenten)}
      \end{array}
      \\[\jot]
    \cup\ \setof{\symbc_\alpha\mathbin{\symb{\in}}\symb{d}}{\alpha<\kappa}\\
    \cup\ \ssetof{\symb{d}\mathbin{\symb{\in}}\symbc_\kappa}\qquad\qquad\qquad\qquad
    \qquad\qquad\qquad\qquad\qquad
    \ifarxived
    \smash{\raisebox{1.44ex}{$\left.\rule[0pt]{0pt}{2.88ex}\right\}$\quad  \assertof{**}}}
    \else
    \smash{\raisebox{0.8zh}{$\left.\rule[0pt]{0pt}{1.6zh}\right\}$\quad  \assertof{**}}}
    \fi
  \end{array}
  $
\end{xitemize}
where $\varphi^*$ is the $\calL^\II$-sentence in \xitemof{x-w-ext-3-0}.

$T$ is a $\calL^\II$-theory of size $\kappa$. 
Similarly to the proof of \Thmof{P-w-ext-0} we can show that $T$ 
is $\LT\kappa$-satisfiable. Thus, by definition of $\kappa$,  $T$ has a model $\gmA$. 
The Mostowski collapse $\gmA^*$  
of the 
model $\gmA$ \wrt\ $\symb{\in}^\gmA$ 
is $V_{\overline{\theta}}$ for some $\overline{\theta}>\kappa$ and 
$\mapping{j}{M}{V_{\overline{\theta}}}$; $a\mapsto (\symbc_a)^{\gmA^*}$
is 
an elementary embedding of $M$ into $V_{\overline{\theta}}$. Because of missing 
first line  of \xitemof{x-w-ext-4}\,\assertof{*} in 
$T$ of \xitemof{x-w-ext-4}, our $T$ in \xitemof{x-w-ext-5} cannot 
guarantee that 
$\kappa$ is the critical point of the elementary embedding but its critical point must be some cardinal
$\leq\kappa$ by \xitemof{x-w-ext-5}\,\assertof{**}, and \xitemof{x-w-ext-4-0} holds. 

This shows that $\kappa$ is weakly sub-extendible. \qedofThm
\qedskip

To place the weakly extendible cardinal in the hierarchy (or a zoo?) of (small) large cardinals, let us 
recall some notions of large cardinals we are going to mention in \Thmof{P-w-ext-1}: 

A cardinal $\kappa$ is said to be {\It strongly unfoldable} if $\kappa=2^{\LT\kappa}$ and, for 
any ordinal $\lambda>\kappa$ and any transitive model $M$ of $\ZFC^-$ \st\ $\kappa\in M$,
$\fnsp{\kappa\GT}{M}\subseteq M$ and $\cardof{M}=\kappa$, there is a transitive 
$N\supseteq V_\lambda$ with $\Elembed{j}{M}{N}{\kappa}$ and $j(\kappa)>\lambda$. 
Here, $\ZFC^-$ denotes the axiom system \ZFC\ without the Power Set Axiom. 

The notion of strongly unfoldable cardinal was introduced by A.\,Villaveces in \cite{andres}. Recently, 
P.\ L\"ucke \cite{luecke2} proved that the strong unfoldability is equivalent to the 
shrewdness, a natural strengthening of the total indescribability which was introduced by 
M.\,Rathjen in \cite{rathjen}.  

Similarly to the argument after the definition \xitemof{x-w-ext-0}, \xitemof{x-w-ext-1} of 
weakly extendible cardinals, 
we can also easily prove that a strongly unfoldable cardinal is weakly compact and hence 
inaccessible in particular.

An uncountable cardinal $\kappa$ is said to be {\It subtle} if, for any club
$C\subseteq\kappa$, and any sequence $\seqof{A_\alpha}{\alpha\in C}$ \st\ 
$A_\alpha\subseteq\alpha$ for all $\alpha\in C$, there are $\alpha$, $\beta\in C$ with
$\alpha<\beta$ \st\ $A_\alpha=A_\beta\cap\alpha$. The notion of subtle cardinal was first 
considered by Jensen and Kunen in \cite{jensen-kunen}. Baumgartner \cite{baum1}, 
\cite{baum2} studied further its combinatorial properties and associated ideals. 

Note that subtle cardinals are compatible with $\uniV=\uniL$: if $\kappa$ is subtle, then so it is in $\uniL$. 

Strongly uplifting cardinals are introduced in Hamkins and Johnstone 
\cite{hamkins-johnstone}: an inaccessible cardinal $\kappa$ is {\It strongly uplifting} if, for 
every $A\subseteq\kappa$ there are arbitrarily large regular $\theta>\kappa$ \st\ 
$\pairof{V_\kappa,\in,A}\prec\pairof{V_\theta,\in, \overline{A}}$ for some
$\overline{A}\subseteq V_\theta$. 
\ifextended{\extendedcolor
\smallskip

The following \LemmaAof{P-w-ext-0-1} and \PropAof{P-w-ext-0-2} were originally used in a direct proof 
of \Thmof{P-w-ext-1},\,\assertof{4}. 

\begin{LemmaA}
  \Label{P-w-ext-0-1} Suppose that $\kappa$ is a subtle cardinal. Then: \wassertof{1}\ \
  $\kappa$ is inaccessible.\smallskip

  \wassert{2} For any club $C\subseteq\kappa$ and any sequence
  $\seqof{A_\alpha}{\alpha\in C}$ \st\ $A_\alpha\subseteq \alpha$ for all $\alpha\in C$, 
  there are inaccessible $\alpha$, $\beta\in C$ with $\alpha<\beta$ and
  $A_\alpha=A_\beta\cap\alpha$. 
\end{LemmaA}\oldmemo{\mbox{L\"oSkoT - annotated.pdf} p.78}
\prf \assertof{1}: $\kappa$ is regular: Suppose that $\kappa$ is a singular cardinal. Say,
$\kappa=\sup_{\xi<\mu}\kappa_\xi$  
for some $\mu<\kappa$ and a continuously and strictly increasing sequence of ordinals \st\
$\kappa_0>\mu$. Let $C=\setof{\kappa_\xi}{\xi<\mu}$ and, for each $\kappa_\xi\in C$, let
$A_{\kappa_\xi}=\ssetof{\xi}$. Then $\seqof{A_\eta}{\eta\in C}$ is a counter example for 
the subtleness of $\kappa$. \smallskip

$\kappa$ is a strong limit: Suppose that $\kappa$ is not a strong limit. Then there is 
$\mu<\kappa$ \st\ $2^\mu\geq\kappa$. Let $C=\kappa\setminus\mu$ and let
$\mapping{f}{C}{\psof{\mu}}$ be an injection. Then $\seqof{f(\xi)}{\xi\in C}$ is a 
counterexample for the subtleness of $\kappa$.\smallskip

\assertof{2}: Let $C\subseteq\kappa$ be a club and let $\seqof{A_\alpha}{\alpha\in C}$ be 
\st\ $A_\alpha\subseteq\alpha$ for all $\alpha\in C$. We want to show that there are 
inaccessible $\alpha$, $\beta\in C$ with $\alpha<\beta$ and 
$A_\alpha=A_\beta\cap\alpha$. 

Since $\kappa$ is inaccessible by \assertof{1}, we may assume that 
\begin{xitemize}
\xitemA[x-w-ext-8-0] 
  all $\alpha\in C$ are strong limit cardinals.
\end{xitemize}
For each $\alpha$,
\begin{xitemize}
\item[] $E_\alpha:=$ even ordinals $\subseteq \alpha$, \quad and\quad 
  $O_\alpha:=$ odd ordinals $\subseteq\alpha$. 
\end{xitemize}
Let $\mapping{f_\alpha}{\alpha}{O_\alpha}$ be the unique order preserving bijection for 
each $\alpha\in C$. Note that 
\begin{xitemize}
\xitemA[x-w-ext-8-1] 
  for all $\alpha$, $\beta\in  C$ with $\alpha<\beta$, we have
  $f_\alpha=f_\beta\restr\alpha$.
\end{xitemize}
Let
\begin{xitemize}
\xitemA[x-w-ext-8-2] 
  $B'_\alpha:=f_\alpha\imageof A_\alpha$, 
\xitemA[x-w-ext-8-3]{}\ifarxived\\[-1.8ex]\else\\[-1zh]\fi
  $B''_\alpha=\left\{\,
  \begin{array}{@{}l}
    \mbox{a cofinal subset }D\mbox{ of }E_\alpha\mbox{ of order-type }\cf(\alpha)
    \mbox{, \st\ }0\mbox{ is the}\\
    \mbox{least element of }D\mbox{ and }\cf(\alpha)\mbox{ is the one but least element of }D,\\
    \qquad\qquad\qquad\mbox{if }\alpha\mbox{ is singular};\\[\jot]
    \rlap{$\emptyset$,}\qquad\qquad\qquad\mbox{otherwise}
  \end{array}
  \right.$\\and
\xitemA[x-w-ext-8-4] $B_\alpha=B'_\alpha\dotcup B''_\alpha$
\end{xitemize}
for all $\alpha\in C$.

Applying the subtleness of $\kappa$ to the sequence $\seqof{B_\alpha}{\alpha\in C}$, let
$\alpha$, $\beta\in C$ be \st\ $B_\alpha=B_\beta\cap\alpha$. By \xitemof{x-w-ext-8-2} (and 
\xitemof{x-w-ext-8-1}), it follows that $A_\alpha=A_\beta\cap\alpha$.

Both of $\alpha$ and $\beta$ must be regular by \xitemof{x-w-ext-8-3} and hence 
inaccessible by \xitemof{x-w-ext-8-0}. 
\qedofLemma

\begin{PropA}
  \Label{P-w-ext-0-2} A cardinal $\kappa$ is subtle if and only if
  \begin{xitemize}
  \xitemA[x-w-ext-9] for any club $C\subseteq\kappa$ and any sequence $\seqof{B_\alpha}{\alpha<\kappa}$ 
    \st\ $B_\alpha\subseteq V_\alpha$ for all $\alpha\in C$, there are $\alpha$, $\beta\in C$ \st\
    $\alpha<\beta$, $\alpha$ and $\beta$ are inaccessible, and
    $B_\alpha=B_\beta\cap V_\alpha$. 
  \end{xitemize}
\end{PropA}
\prf It is clear that \xitemof{x-w-ext-9} implies that $\kappa$ is subtle. To prove the 
converse, assume that $\kappa$ is subtle. Let $\mapping{f}{\kappa}{V_\kappa}$ be a 
bijection \st\ $f\imageof\alpha\subseteq V_{\alpha}$ for all $\alpha<\kappa$ and 
$f\imageof{\alpha}=V_\alpha$ if and only if $\cardof{\alpha}=\cardof{V_\alpha}$ (this is 
possible by \Lemmaof{P-w-ext-0-1},\,\assertof{1}).

Let $A_\alpha:=f^{-1}\imageof B_\alpha$ for $\alpha\in C$. By 
\Lemmaof{P-w-ext-0-1},\,\assertof{2}, there are inaccessible $\alpha$, $\beta<\kappa$ with
$\alpha<\beta$ and $A_\alpha=A_\beta\cap\alpha$. Noticing that the mappings 
$\mapping{f\restr\alpha}{\alpha}{V_\alpha}$ and $\mapping{f\restr\beta}{\beta}{V_\beta}$ 
are bijective, we obtain
\begin{xitemize}
\item[] $B_\alpha=f\imageof{A_\alpha}=f\imageof{(A_\beta\cap \alpha)}
  =(f\imageof{A_\beta})\cap(f\imageof\alpha)=B_\beta\cap V_\alpha$. \qedofProp
\end{xitemize}
}\fi 

\begin{Thm}\Label{P-w-ext-0-1-0}
  {\rm (Hamkins and Johnstone \cite{hamkins-johnstone} Theorem 7) }\quad Any subtle cardinal is a 
    stationary limit of strongly uplifting cardinals. \qed
\end{Thm}

\begin{Thm}
  \Label{P-w-ext-1} \wassertof{1}\ \ 
  If $\kappa$ is weakly extendible, then there is a weakly compact
  $\lambda>\kappa$. On the other hand, strong unfoldability of $\kappa$ does not imply the 
  existence of inaccessible $\lambda>\kappa$. \smallskip

  \wassert{2} If $\kappa$ is weakly extendible and $\nu$ is the first inaccessible 
  cardinal above $\kappa$ (which exists by {\rm\assertof{1}}), then
  $V_\nu\modelof{\kappa\mbox{ is strongly unfoldable but not weakly extendible}}$. Also
  $V_\nu\modelof{\mbox{there is no inaccessible cardinal above }\kappa}$. It follows 
  that $\ZFC$ $+$ ``\,there is a weakly extendible cardinal\/'' proves 
  $consis(\nnumeralof{\ZFC} +\numeralof{\xmbox{\,there is a strongly unfoldable cardinal\,}})$.\smallskip

  \wassert{3} {\rm (Boney, Dimopoulos, Gitman, and Magidor \cite{boneyetal}, Proposition 
    4.8)}\quad If $\kappa$ is a strongly uplifting cardinal then $\kappa$ is weakly 
  extendible and $\kappa$ is a stationary 
  limit of weakly extendible cardinals. \smallskip

  \wassert{4} If $\kappa$ is a subtle cardinal, then $\kappa$ is a stationary limit of 
  weakly  extendible cardinals. \memo{A model in which $\kappa$ is subtle but not extendible?}
\end{Thm}
\prf \assertof{1}: Suppose that $\kappa$ is a weakly extendible cardinal. By 
\Lemmaof{P-w-ext-a},\,\assertof{1}, $\kappa$ is weakly compact.

Let
$M\prec V_\theta$ for a $\theta>\kappa$ be \st\
$V_\theta\modelof{\kappa\mbox{ is weakly compact}}$, $\kappa+1\subseteq M$ and
$\cardof{M}=\kappa$. Let $\overline{\theta}>\kappa$ and $j$ be \st\
$\EElembed{j}{M}{V_{\overline{\theta}}}{\kappa}$. 

By elementarity, $V_{\overline{\theta}}\modelof{j(\kappa)\mbox{ is weakly compact}}$. Hence 
$j(\kappa)>\kappa$ is really weakly compact.

The second assertion follows from \assertof{2}. \smallskip

\assertof{2}: Suppose that $\kappa$ is weakly extendible and $\nu>\kappa$ is the first 
inaccessible cardinal above $\kappa$.

Since $\kappa=2^{\LT\kappa}$, we also have $V_\nu\models\kappa=2^{\LT\kappa}$. 

Suppose that $M\in V_\nu$ is a transitive model of $\ZFC^-$ \st\ $\kappa\in M$,
$\fnsp{\kappa\GT}{M}\subseteq M$ and $\cardof{M}=\kappa$. Let $\kappa<\lambda<\nu$
($=\On^{V_\nu}$). We have to show that, in $V_\nu$, there is an elementary embedding of $M$ 
into a target model satisfying the conditions in the definition of strong unfoldability.

Let $\kappa<\theta<\nu$ be \st\ $M\in V_\theta$ and let $M^*\prec V_\theta$ be \st\ $M\in M^*$,
$M\subseteq M^*$, $\kappa+1\subseteq M^*$, and $\cardof{M^*}=\kappa$. By weak extendibility 
of $\kappa$, there is $\overline{\theta}>\kappa$ and $j$ \st\
$\EElembed{j}{M^*}{V_{\overline{\theta}}}{\kappa}$. Note that $j(\kappa)\geq\nu$ since 
$j(\kappa)$ must be inaccessible.

Let $N\prec V_{\overline{\theta}}$ be \st\ $V_\lambda$, $j\imageof M^*\subseteq N$ and
$\cardof{N}<\nu$. Let $\mapping{m}{N}{N^*}$ be the Mostowski collapse. Then we have
$N^*\in V_\nu$, and, letting $j^*:=m\circ j$, \,$\EElembed{j^*}{M^*}{N^*}{\kappa}$. Note that 
we also have $j^*\in V_\nu$. 

Since
$V_\lambda\subseteq N$, we have $m\restr V_\lambda=id_{V_\lambda}$ and hence
$V_\lambda\subseteq N^*$. Also, it follows that $j^*(\kappa)\geq\lambda$.

By elementarity, $N^*\modelof{j^*(M)\mbox{ is transitive}}$. Since $N^*$ is transitive it 
follows that $j^*(M)$ is really transitive. Also, by elementarity
$N^*\models[j^*(M)]^{\LT j^*(\kappa)}\subseteq j^*(M)$. It follows that
$V_\lambda\subseteq j^*(M)$.

It is easy to check that $\Elembed{j^*\restr M}{M}{j^*(M)}{\kappa}$ and all of these 
objects are elements of $V_\nu$.

This shows that $\kappa$ is strongly unfoldable in $V_\nu$. 

Clearly, $V_\nu\modelof{\mbox{there is no inaccessible cardinal }\GT\kappa}$ and this shows 
the second half of \assertof{1}.

Since $V_\nu$ as above is also a model of \ZFC, we obtain
$consis(\nnumeralof{\ZFC} +\numeralof{\xmbox{\,there is a strongly unfoldable cardinal\,}})$. 
\smallskip

\assertof{3}: For the proof of the first claim see \cite{boneyetal} (noticing 
\Thmof{P-w-ext-0}). The second claim is obtained by a 
slight modification of the proof given in \cite{boneyetal}: Suppose that $\kappa$ is a 
strongly uplifting cardinal and $C$ is a club subset of $\kappa$. By the first claim, 
$\kappa=\wcn(\calL^\II_{\kappa,\omega})$. Let $\theta$ be sufficiently large \st\
$V_\theta\modelof{\kappa=\wcn(\calL^\II_{\kappa,\omega})}$ and there is $\overline{C}\subseteq\theta$ \st\
$\pairof{V_\kappa,\in,C}\prec\pairof{V_\theta,\in,\overline{C}}$.

Since
$\pairof{V_\theta,\in,\overline{C}}
\modelof{\symb{C}(\kappa)\mbox{ and }\kappa=\wcn(\calL^{\II}_{\kappa,\omega})}$, we have 
\begin{xitemize}
\item[] 
  $\pairof{V_\theta,\in,\overline{C}}\modelof{\xmbox{there is a cardinal }\delta\mbox{ in }
    \symb{C}\mbox{ with }\delta=\wcn(\calL^{\II}_{\delta,\omega})}$
\end{xitemize}
(we denote with $\symb{C}$ 
the unary predicate symbol corresponding to $C$ and $\overline{C}$), it follows that
\begin{xitemize}
\item[] 
  $\pairof{V_\kappa,\in,C}\modelof{\xmbox{there is a cardinal }\delta\mbox{ in }
    \symb{C}\mbox{ with }\delta=\wcn(\calL^{\II}_{\delta,\omega})}$ 
\end{xitemize}
by elementarity. Since $\kappa$ is strongly uplifting, it follows 
that ``$\delta=\wcn(\calL^{\II}_{\delta,\omega})$'' holds in $\uniV$ for $\delta$ as above. 
\smallskip

\assertof{4}: follows from \Thmof{P-w-ext-0-1-0} and \assertof{3}. 
\ifextended{\extendedcolor Suppose that $\kappa$ is a subtle cardinal 
  and $C\subseteq\kappa$ is a club. By \Thmof{P-w-ext-0-1-0}, there is $\lambda\in\Lim(C)$ 
  \st\ $\lambda$ is strongly uplifting. Now since $C\cap\lambda$ is a club in $\lambda$, 
  there is a weakly extendible $\mu\in C\cap\lambda$ by \assertof{3}. 
}\fi 
\qedofThm


\Thmof{P-w-ext-0-0}  implies that $\wcn(\calL^\II)$ cannot be a small cardinal. 
This follows from the next observation which should be well-known:

\begin{Lemma}
  \Label{P-w-ext-3} Suppose that \ifextended{\extendedcolor\ \ \ifarxived\daimaru{1}\else \fi} 
  \fi 
  $M\prec V_\theta$, $\kappa+1\subseteq M$, $N$ and \ifextended{\extendedcolor\ \ \ifarxived\daimaru{2}\else \fi}\fi
  $\EElembed{j}{M}{N}{\mu}$ for $\mu\leq\kappa$. Then, $\mu$ is \ifextended{\extendedcolor\wassertof{1} a regular 
  cardinal, \wassertof{2} }\fi weakly inaccessible, and 
  \ifextended{\extendedcolor\\ \wassertof{3} }\else further, \fi weakly Mahlo, weakly hyper 
  Mahlo, etc.\ifextended\else\\ \qed\fi
\end{Lemma}
\ifextended{\extendedcolor \prf
\Wolog, we may assume that $N$ is transitive.

\assertof{1}: Suppose toward a contradiction 
that $\mu$ is singular. Then, by the elementarity \ifarxived\daimaru{1}\else \fi, there are $\nu<\mu$ and $f\in M$ \st\
$M\modelof{\mapping{f}{\nu}{\mu}\mbox{ is cofinal}}$. Since
$j(f)=j(f)\restr j(\nu)=j(f)\restr\nu=f$ it follows, by the elementarity 
\ifarxived\daimaru{2}\else \fi, that $N\modelof{\mapping{f}{\nu}{j(\mu)}\mbox{ is cofinal}}$. Thus 
$\mu=j(\mu)$. This is a contradiction. \smallskip

\assertof{2}: By \assertof{1} (and since $\mu$ cannot be $\omega$), it is enough to show 
that $\mu$ is not a successor cardinal. Suppose it were, say $\mu=\nu^+$. $j(\mu)=\mu$ 
by definition of $\mu>\nu$.  By the elementarity \ifarxived\daimaru{1}\else \fi, we have
\begin{xitemize}
\item[] 
  $M\modelof{\forall x<\mu\,(x\not=0\ \rightarrow\ \exists y(\mapping{y}{\nu}{x}\mbox{ is a surjection}))}$. 
\end{xitemize}

By the elementarity \ifarxived\daimaru{2}\else \fi, it follows that
\begin{xitemize}
\item[] 
  $N\modelof{\forall x<j(\mu)\,(x\not=0\ \rightarrow\ \exists y(\mapping{y}{\nu}{x}\mbox{ is a surjection}))}$. 
\end{xitemize}
Thus we have $\mu\leq j(\mu)\leq\nu^+=\mu$ and hence $j(\mu)=\mu$. This is a 
contradiction to the choice of $\mu$.\smallskip

\assertof{3}: Suppose that $C\in M$ is a club subset of $\mu$. Then $j(C)\cap \mu= C$ and
$N\modelof{j(C)\mbox{ is a club subset of }j(\mu)}$. This and \assertof{2} together with the 
elementarity \ifarxived\daimaru{2} \else   \fi imply 
$N\modelof{\mu\in j(C)\mbox{ and }\mu\mbox{ is a weakly inaccessible cardinal}}$. 
Again by the elementarity \ifarxived\daimaru{2}\else  \fi and since $C$ was arbitrary, we obtain
\begin{xitemize}
\item[] $
  \begin{array}{r@{}l}
    M\modelof{\forall_x\,({}&x\mbox{ is a club subset of }\mu\ \rightarrow\\ 
      &x\mbox{ contains an element which is weakly inaccessible})}.
  \end{array}
  $
\end{xitemize}

By the elementarity \ifarxived\daimaru{2}\else \fi, the same statement holds in $V_\theta$ and hence also in $V$. This 
shows that $\mu$ is weakly Mahlo. The same argument can be repeated to show that $\mu$ is 
weakly hyper Mahlo, weakly hyper hyper Mahlo, etc. 
\qedofLemma}\fi 

\begin{Prop}
  \Label{P-w-ext-4}
  \wassertof{1} If $\wcn(\calL^\II)<\infty$, $\wcn(\calL^\II)$ is greater than 
  the first weakly Mahlo cardinal, first weakly hyper Mahlo cardinal, etc. \smallskip

  \wassert{2} $\continuum\notin\WCS(\calL^\II)$, $(\continuum)^+\not\in\WCS(\calL^\II)$, etc. 
\end{Prop}
\prf \assertof{1}: By \Thmof{P-w-ext-0-0} and \Lemmaof{P-w-ext-3}.

\assertof{2}: $T:={}
\begin{array}[t]{@{}l}
  \mbox{``the }\calL^\II\mbox{-theory of }\pairof{\psof{\omega}, a, \in}_{a\in\psof{\omega}}\mbox{''}
  \ \cup\ \setof{\symb{d}\not=\symb{c}_a}{a\in\psof{\omega}}
\end{array}
$\\
is a counter-example to $\continuum\in\WCS(\calL^\II)$. 

For $n\in\omega\setminus 1$, let $\mapping{f}{\psof{\omega}}{\continuum}$ be a bijection, 
and let 
\begin{xitemize}
\item[] 
  $\calF_n:=\setof{g}{\mapping{g}{(\continuum)^{+k}}{\alpha}
  \mbox{ is a surjection for }k\in n,\,0<\alpha<(\continuum)^{+(k+1)}}$. 
\end{xitemize}
Then 
\begin{xitemize}
\item[] 
  $T_n:={}
\begin{array}[t]{@{}l}
  \mbox{``the }\calL^\II\mbox{-theory of }\\
  \qquad\pairof{\psof{\omega}\cup (\continuum)^{+n}, 
    \psof{\omega}, \continuum, (\continuum)^{+1}\ctentenc (\continuum)^{+(n-1)}, f, g, \in}_{g\in\calF_n}\mbox{''}\\
  \ \cup\ \setof{\symb{d}\not=\symb{c}_a}{a\in\psof{\omega}\cup (\continuum)^{+n}}
\end{array}$
\end{xitemize}
is a counter example to $(\continuum)^{+n}\in\WCS(\calL^\II)$. 
\qedofProp
\section{General characterization of weakly $\calL$-compact cardinals}
\Label{comp-card}
We show in this section that \Thmof{P-w-ext-0} can be generalized to a wider class of logics.

In the following, $\ZC^-$ denotes the Zermelo set-theory with the Axiom of Choice minus the 
Power-set Axiom.

Suppose that $N\models\ZC^-$. 
For a signature $\calS\in N$ with $\calS=\pairof{\symbc_i, \symb{f}_j, \symb{r}_k}_{i\in \kappa_0,j\in\kappa_1,k\in\kappa_2}$ 
where $\symbc_i$, is a constant symbol, $\symb{f}_j$ an $m_j$-place function symbol, and $\symb{r}_k$ an $n_k$-place 
relation symbol for each $i\in\kappa_0$, $j\in\kappa_1$, $k\in\kappa_2$; 
$\kappa_0+1$, $\kappa_1+1$, $\kappa_2+1\subseteq N$, and 
\begin{xitemize}
\xitem[x-LS-24-a] 
  $N\modelof{\gmA\xmbox{ is a structure in the signature }
  \ssetof{\smash{\symb{c}}_i,\smash{\symb{f}}_j, 
    \smash{\symb{r}}_k}_{i\in\kappa_0,j\in\kappa_1,k\in\kappa_2}}$ 
\end{xitemize}
we denote with 
$\gmA^N$ the structure 
\begin{xitemize}
\xitem[x-LS-24-0] ${\darkred\gmA^N}:=\pairof{A\cap N, c_i, f_j\restr (A\cap N)^{m_j},
  r_k\cap (A\cap N)^{n_k}}_{i\in \kappa_0,j\in\kappa_1,k\in\kappa_2}$ 
\end{xitemize}
where $c_i=(\symbc_i)^\gmA$, $f_j=(\smash{\symb{f}}_j)^\gmA$, etc. 

Note that, 
\begin{xitemize}
\xitem[x-LS-24] 
  if $\isom{m}{N}{N_0}$ is the Mostowski collapse, then $\gmA^N\cong m(\gmA)=(m(\gmA))^{N_0}$. 
\end{xitemize}

For a logic $\calL$, if $N$ is \st\ $N\models\ZC^-$ and $N$ contains all 
parameters needed to define $\calL$, we shall say that 
$N$ is {\It$\calL$-truthful\/} if, for all structures $\gmA$ as 
above (in connection with this $N$), 
$N\modelof{\gmA\models_\calL\varphi}$ is equivalent to 
$\gmA^N\models_{\calL}\varphi$. By \xitemof{x-LS-24}, if $N$ is $\calL$-truthful, then its 
Mostowski collapse $N_0$ is also $\calL$-truthful.

Note that, for a sentence $\varphi$ in a proper logic $\calL$, there is a first-order formula $\varphi_\ast$ \st\
\begin{xitemize}
\xitem[x-LS-24-1] 
  $N\modelof{\gmA\models_\calL\varphi}$\ \ $\Leftrightarrow$\ \ $N\models\varphi_\ast(\gmA)$,
\end{xitemize}
and an $\calL$-formula $\varphi_{\ast\ast}$ \st\ 
\begin{xitemize}
\xitem[x-LS-24-2] 
  $\gmA^N\models_\calL\varphi$\ \ $\Leftrightarrow$\ \ 
  $N\models_\calL\varphi_{\ast\ast}(\gmA)$.
\end{xitemize}

Thus, letting
\begin{xitemize}
\xitem[x-LS-24-3] 
  $\Phi^*_\calL:=\setof{\forall_x(\mbox{``}x\mbox{ is a structure in the signature of $\varphi$''}\rightarrow
  (\varphi_\ast(x)\leftrightarrow\varphi_{\ast\ast}(x)))\\\mbox{}\hfill}{\varphi\mbox{ is 
        an $\calL$-formula}}$,\ \mbox{}
\end{xitemize}
we have
\begin{xitemize}
\xitem[x-LS-24-4] 
  $N\models\Phi^*_\calL$\ \ $\Leftrightarrow$\ \  $N$ is $\calL$-truthful\quad
\end{xitemize}
for any $N$ with $N\models\ZC^-$

Let us call here a logic $\calL$ {\It finitary} if the set of free variables in any $\calL$-formula 
is finite, the set of all $\calL$-formulas of given signature $S$ of cardinality $\leq\kappa$, 
for an infinite $\kappa$ has size $\leq\kappa$, 
and, for any infinite ordinal $\theta$ and $\varphi\in V_\theta$,
$V_\theta\modelof{\varphi\mbox{ is an }\calL\mbox{-formula}}$ if and only if 
$\varphi$ is (really) an $\calL$-formula.

$\calL_{\kappa,\omega}$ for an uncountable 
cardinal $\kappa$ is an example of non finitary logic (since the size of the set of 
formulas in a signature can exceed the size of the signature).

\begin{Thm}\Label{P-comp-2} \wassertof{1}\ \  Suppose that $\calL$ is a finitary proper logic 
  \st\ 
  \begin{xitemize}
  \xitem[x-comp-6] $V_\theta$ for all regular uncountable $\theta$ is $\calL$-truthful; and
  \xitem[x-comp-7] ``$\symbin$ is well-founded'' is 
    expressible by a formula $\varphi^*_\calL$ in $\calL$. 
  \end{xitemize}

  Then a \memo{regular?} cardinal $\kappa$ is weakly $\calL$-compact (i.e.\
  $\kappa\in\WCS(\calL)$\/)\ \ $\Leftrightarrow$ 
  \begin{xitemize}
  \xitem[x-comp-8] for any regular $\theta\geq\kappa$ and $M\prec V_\theta$ \st\
    $\kappa+1\subseteq M$, $\cardof{M}=\kappa$, there are $j$, $N$ \st\
    $\kappa+1\subseteq N$, $\eelembed{j}{M}{N}$,
    $j(\kappa)>\min(\On^N\setminus\sup(j\imageof\kappa))$, and $N$ 
    is $\calL$-truthful.
  \end{xitemize}

  \wassert{2} Suppose that $\calL^*$ is a logic obtained from a finitary proper logic $\calL$
  which satisfies \xitemof{x-comp-6} and \xitemof{x-comp-7}, by extending $\calL$ by taking the closure 
  of the set of $\calL$ formulas \wrt\ infinitary conjunction and disjunction of set of formulas of 
  size $\LT\kappa$ and first order logical operations.
  
  Then $\kappa$ is weakly $\calL^*$-compact
  \ \ $\Leftrightarrow$\ \ $2^{<\kappa}=\kappa$ and 
  \begin{xitemize}
  \xitem[x-comp-8-0] for any regular $\theta\geq\kappa$ and $M\prec V_\theta$ \st\
    $\kappa+1\subseteq M$, $\cardof{M}=\kappa$, there are $j$, $N$ \st\
    $\EElembed{j}{M}{N}{\kappa}$,\/\footnotemark
    and
    $N$ is $\calL$-truthful.  
  \end{xitemize}
\end{Thm}
\footnotetext{For the notation $\EElembed{j}{M}{N}{\kappa}$ see directly after \xitemof{x-w-ext-1}.}
\prf \assertof{1}: ``$\Leftarrow$'': Assume that \xitemof{x-comp-8} holds for $\kappa$. Let 
$T$ be a $\LT\kappa$-satisfiable $\calL$-theory of signature of size $\leq\kappa$. 
Since $\calL$ is finitary, we may assume that $\cardof{T}=\kappa$. \Wolog, we may further assume that 
$T$ is coded as a subset of $\kappa$. 
Let $\theta$ be 
sufficiently large regular cardinal \st\
$V_\theta\modelof{T\mbox{ is }\LT\kappa\mbox{-satisfiable}}$.

Let $M\prec V_\theta$ be \st\ $T\in M$, $\kappa+1\subseteq M$ (note that this implies
$T\subseteq M$), and $\cardof{M}=\kappa$. By \xitemof{x-comp-8}, there are $j$ and $N$ \st\  
\begin{xitemize}
\xitem[x-comp-9] 
  $\kappa+1\subseteq N$,
\xitem[x-comp-10] $\eelembed{j}{M}{N}$,
\xitem[x-comp-11] $j(\kappa)>\min(\On^N\setminus\sup(j\imageof\kappa))$, and
  \xitem[x-comp-12] $N$ is $\calL$-truthful.
\end{xitemize}
By elementarity \xitemof{x-comp-10}, we have
\begin{xitemize}
\xitem[x-comp-13] 
  $N\modelof{j(T)\mbox{ is a }\LT j(\kappa)\mbox{-satisfiable }\calL\mbox{-theory}}$.  
\end{xitemize}
Let $\alpha^*:=\min(\On^N\setminus\sup(j\imageof\kappa))$ and $T^*:=j(T)\cap \alpha^*$. 
Then $T^*\in N$ and $j\imageof T\subseteq T^*$. By \xitemof{x-comp-11} and \xitemof{x-comp-13}, it follows 
that $N\modelof{\gmA\models T^*}$ for some structure $\gmA\in N$. By \xitemof{x-comp-12}, it follows that
$\gmA^N\models T^*$ and hence $\gmA^N\models j\imageof T$. Thus, by renaming the components 
of the structure $\gmA^N$, we obtain a model of $T$. \smallskip
\memo{This argument works by what is called substitution property in \cite{MTL} p.30 under  
  the assumption the coding $T\subseteq\kappa$ is done appropriately.}

``$\Rightarrow$'': Assume that $\kappa$ is weakly $\calL$-compact. 

Suppose that $\theta\geq\kappa$ is a regular cardinal, and $M$ is \st\ $M\prec V_\theta$,  
    $\kappa+1\subseteq M$, and $\cardof{M}=\kappa$.
Let $T$ be the $\calL$-theory defined by:
\begin{xitemize}
\xitem[x-comp-14] ${\darkred T}:={}
  \begin{array}[t]{@{}l}
    \ssetof{\varphi^*_\calL}
    \\ 
    \cup\ \setof{\varphi(\symbc_{a_0}\ctenten)}{{}
      \begin{array}[t]{@{}l}
        \varphi
        \mbox{ is a first-order formula in the signature }\ssetof{\symb{\in}},\\
        a_0\ctenten\in M\mbox{ and }M\models\varphi(a_0\ctenten)}
      \end{array}
      \\[\jot]
    \cup\ \setof{\symbc_\alpha\mathbin{\symb{\in}}\symb{d}}{\alpha<\kappa}\\
    \cup\ \ssetof{\symb{d}\mathbin{\symb{\in}}\symbc_\kappa}\\
    \cup\ \Phi^*_\calL
  \end{array}
  $
\end{xitemize}
where $\varphi^*_\calL$ is an $\calL$-formula as in \xitemof{x-comp-7} and $\Phi^*_\calL$ a set 
of $\calL$-formulas in \xitemof{x-LS-24-3}. 

$T$ is $\LT\kappa$-satisfiable: by \xitemof{x-comp-6} and \xitemof{x-LS-24-4}, $V_\theta$ can be expanded to a 
model $\gmA_0$ of any subset $T_0$ of $T$ of size $\LT\kappa$ by letting $(\symbc_a)^\gmA:=a$ for 
all constant  symbols of this form appearing $T_0$, and
$(\symb{d})^\gmA:=\sup(\setof{\alpha<\kappa}{\symbc_\alpha\mbox{ appears in }T_0})$. Note 
that $\sup(\setof{\alpha<\kappa}{\symbc_\alpha\mbox{ appears in }T_0})<\kappa$ since 
$\kappa$ is regular.

By assumption, there is a model $\gmB$ of $T$. By $\varphi^*_\calL\in T$, 
we can take the Mostowski collapse $\gmB^*$ of $\gmB$ (\wrt\ $\symb{\in}^\gmB$). 
We then have $\symbin^{\gmB^*}={\in}$ and $\ulsetof{\gmB^*}$ is a transitive set. 

By \xitemof{x-LS-24-4} and since $\Phi^*_\calL\subseteq T$ , 
$N:=\ulsetof{\gmB^*}$ is $\calL$-truthful. 

Let $\mapping{j}{M}{N}$ be defined by $j(a):=(\symbc_a)^{\gmB^*}$ for $a\in M$. Then these $N$, $j$ are 
as desired in \xitemof{x-comp-8}. \smallskip

\assertof{2}: This can be proved similarly to \assertof{1}. For the direction ``$\Rightarrow$'', note the 
following \Lemmaof{P-comp-2-0}. Also, note that $2^{\LT\kappa}=\kappa$ follows from the 
assumption of the weakly $\calL^*$-compactness of $\kappa$ since $\kappa$ is then at 
least weakly compact. 

For the proof 
of the direction ``$\Leftarrow$'', consider
\begin{xitemize}
\item[] ${\darkred T'}:={}
  \begin{array}[t]{@{}l}
    \ssetof{\varphi^*_\calL}
    \\ 
    \cup\ \setof{\varphi(\symbc_{a_0}\ctenten)}{{}
      \begin{array}[t]{@{}l}
        \varphi
        \mbox{ is a first-order formula in the signature }\ssetof{\symb{\in}},\\
        a_0\ctenten\in M\mbox{ and }M\models\varphi(a_0\ctenten)}
      \end{array}
      \\[\jot]
    \cup\ \setof{\forall_x(x\mathbin{\symb{\in}}\symbc_\alpha\ \leftrightarrow\ 
      \llor_{\beta<\alpha}x\equiv\symbc_\beta)}{\alpha<\kappa}\\
    \cup\ \setof{\symbc_\alpha\mathbin{\symb{\in}}\symb{d}}{\alpha<\kappa}\\
    \cup\ \ssetof{\symb{d}\mathbin{\symb{\in}}\symbc_\kappa}\\
    \cup\ 
    \Phi^*_{\calL^*}
  \end{array}
  $
\end{xitemize}
in place of $T$ in the proof of \assertof{1}. 
\qedofThm

\begin{Lemma}
  \Label{P-comp-2-0} Suppose that $\calL$ and $\calL^*$ are as in \Thmof{P-comp-2}. If $M$ 
  is a model of a sufficiently large fragment of set theory, $M\supseteq\kappa+1$, and $M$ is $\calL$-truthful, 
  then $M$ is also $\calL^*$-truthful. \qed
\end{Lemma}

\ifextended{\extendedcolor
\begin{LemmaA}
  \Label{P-comp-3-0} Suppose that $\calL$ is a finitary proper logic and $M$ is a 
  transitive model of $\ZC^-$ \st\ \  
\ixitem[x-comp-16] $M$ is $\calL$-truthful,  
  $N\in M$, and\ \ \ixitemx[x-comp-17] $M\modelof{N\xmbox{ is an $\calL$-truthful model 
    of }\ZC^-}$.\footnotemarkA\ Then $N$ is $\calL$-truthful. 
\end{LemmaA}
\footnotetextA{In general, this can be formulated by infinitely many conditions: one for 
  each $\calL$-formula.}
\prf Let $\gmA\in N$ and $\varphi$ is a $\calL$-formula in $N$ in the signature of $\gmA$. 
Then we have
\begin{xitemize}
\item[] $N\modelof{\gmA\models_\calL\varphi}$\
  ${\obecause{\Leftrightarrow}{}{by \xitemof{x-comp-16} and since $\calL$ contains the 
      first-order logic}}$\ 
  $M\modelof{N\modelof{\gmA\models_\calL\varphi}}$\ 
  $\smash{\ubecause{\Leftrightarrow}{4pt}{by \xitemof{x-comp-17}}}$\  
  $M\modelof{\gmA^N\models_\calL\varphi}$\\
  ${\ubecause{\Leftrightarrow}{-2pt}{by \xitemof{x-comp-16}}}$\ 
  $\gmA^N\models_\calL\varphi$.\ifarxived\\[-2.7ex]\else\\[-1.5zh]\fi\qedofLemma
\end{xitemize}}\fi

\Thmof{P-comp-2} is a generalization of \Thmof{P-w-ext-0}. This can be seen in the following:  
\begin{Lemma}
  \Label{P-comp-4} If a transitive set $N$ (seen as an $\in$-structure) is elementary 
  equivalent (in first-order logic) to $V_\theta$ for  
  some limit ordinal $\theta$ and $N$ is $\calL^\II$-truthful then $N=V_{\theta'}$ for some 
  limit ordinal $\theta'$. 
\end{Lemma}
\prf Note that, by elementarity, 
\begin{xitemize}
\item[] $N\modelof{\forall_x\exists_\alpha\exists_f\,(\,
  \begin{array}[t]{@{}l}
    \mbox{``}\alpha\mbox{ is an ordinal''}\ \land\ \mbox{``}f\mbox{ is a mapping on }\alpha\mbox{''}\\
    \land\ \mbox{``}f\mbox{ is an initial segment of the cumulative hierarchy''}\\
    \land\ \mbox{``}x\in f(\beta)\mbox{ for some }\beta<\alpha)}. 
  \end{array}$
\end{xitemize}

Thus, it is enough to show that, if $A\in N$, then $\psof{A}\subseteq N$. 
Let $\gmA\in N$ be \st\ 
$N\modelof{\gmA=\pairof{\psof{A}\cup A, A,\in}}$. Since $N$ is transitive, we have 
\ixitem[x-comp-17-0] $\gmA^N=\gmA$. 

\begin{xitemize}
\item[] $N\modelof{\gmA\models_{\calL^\II}\forall X(\forall_y(y\varin X\rightarrow y\symbin\symb{A})
  \ \rightarrow\ \exists_z\forall_y\,(y\symbin z\leftrightarrow y\varin X))}$, 
\end{xitemize}
since $N$ is elementary equivalent to $V_\theta$. 

It follows that 
\begin{xitemize}
\item[] $\gmA\models_{\calL^\II}\forall X(\forall_y(y\varin X\rightarrow y\symbin\symb{A})
  \ \rightarrow\ \exists_z\forall_y\,(y\symbin z\leftrightarrow y\varin X))$ 
\end{xitemize}
by $\calL^\II$-truthfulness of $N$ (and \xitemof{x-comp-17-0}).

It follows that $\ulsetof{\gmA^N}=\psof{A}\cup A$. Thus, 
$\psof{A}\subseteq N$. \qedofLemma
\smallskip

$\wcn(\calL)$ for a finitary logic $\calL$ as in \Thmof{P-comp-2},\,\assertof{1} is always fairly 
large since the following analog of \Propof{P-w-ext-4} holds for such a logic:

\begin{Prop}
  \Label{P-comp-4-0} Suppose that $\calL$ is a finitary proper logic satisfying 
  \xitemof{x-comp-6} and \xitemof{x-comp-7}. 

  \wassertof{1} If $\wcn(\calL)<\infty$, $\wcn(\calL)$ is greater than 
  the first weakly Mahlo cardinal, first weakly hyper Mahlo cardinal, etc. \smallskip

  \wassert{2} $\continuum\notin\WCS(\calL)$, $(\continuum)^+\not\in\WCS(\calL)$, etc. 
\end{Prop}
\prf \assertof{1}: By \Thmof{P-comp-2},\,\assertof{1} and \Lemmaof{P-w-ext-3}. 

\assertof{2}: The theories similarly to $T$ and $T_n$ in the proof of \Propof{P-w-ext-4}
with $\calL^\II$-theories of the structures replaced by corresponding $\calL$-theories 
show the inequality. \qedofProp






\section{Weak compactness of stationary logic}
\Label{stat}
As is already noted, $\lsn(\calL^\II)$ is the least 
supercompact cardinal by one of Magidor's Theorems (see \Thmof{P-intro-7-0}). In contrast, 
$\lsn(\calL^{\aleph_0,\II})$ is $(\continuum)^+$ while $\lsn(\calL^{\aleph_0}_{\stat})$ can be  
(consistently) $\leq\continuum$ \memo{$\MA^{+\omega_1}(\sigma\mbox{-closed})$ $+$
  $\continuum=\aleph_2$ implies that $\lsn(\calL^{\aleph_0}_{\stat})$ is $\continuum$.\\
  weakly compact ¨ stat reflection Jech \mbox{Millennium Cor.17.20}. \Lemmaof{P-stat-0} 'Ì
  –½'è"ԍ†—v•ÏX}
though the consistency strength of
$\LSS(\calL^{\aleph_0}_\stat)\not=\emptyset$ is still rather high (see the following 
\Lemmaof{P-stat-0},\,\assertof{4} and \assertof{5}).  The size of
$\lsn(\calL^{\aleph_0}_\stat)$ is wildly undecided under ZFC\ (see \Lemmaof{P-stat-0},\,\assertof{10}). 

We shall call a cardinal $\kappa$ {\It $\cdot^{\aleph_0}$-closed}\/ if 
$\mu^{\aleph_0}<\kappa$ holds for all $\mu<\kappa$. 

For a class $\calP$ of \pos, a cardinal $\kappa$ is said to be $\calP$-generically 
supercompact ($\calP$-gen.\ supercompact, for short) if for, any $\lambda\geq\kappa$, there 
is a \po\ $\poP\in\calP$ \st, for a $(\uniV,\poP)$-generic $\genG$, there are classes $j$ 
and $M$ (in $\uniV[\genG]$) \st\ 
\begin{xitemize}
\xitem[x-stat-a] 
  $\uniV[\genG]\modelof{\Elembed{j}{\uniV}{M}{\kappa}}$, 
  $j(\kappa)>\lambda$, and $j\imageof\lambda\in M$
\end{xitemize}
where we are using here the convention in the \footnoteof{fn-0}. 

\begin{Prop}
  \phantomsection
  \Label{P-stat-0} \wassertof{1}\ \
  $\LSS(\calL^{\aleph_0,\II})=\setof{\kappa}{\kappa\mbox{ is }\cdot^{\aleph_0}\mbox{-closed\/}}$.\smallskip

  \wassert{2} $\lsn(\calL^{\aleph_0,\II})=(\continuum)^+$.\smallskip

  \wassert{3} If\/ $\square_\kappa$ holds then $\lsn(\calL^{\aleph_0}_\stat)>\kappa^+$. 
  \smallskip

  \wassert{4} If\/ $\square_\kappa$ holds for class many $\kappa$, then
  $\LSS(\calL^{\aleph_0}_\stat)=\emptyset$. \smallskip

  \wassert{5} 
  $\MA^{+\omega_1}(\sigma\mbox{-closed\/})$ implies
  $\lsn(\calL^{\aleph_0}_\stat)=\aleph_2$. 
  In particular, $\lsn(\calL^{\aleph_0}_\stat)=\continuum$ is consistent (modulo a 
  large cardinal, e.g.\ a supercompact). Note that, by \assertof{2},
  $\lsn(\calL^{\aleph_0,\II}_\stat)\leq\continuum$ is inconsistent. 
  \smallskip 

  \wassert{6} $\lsn(\calL^{\aleph_0}_\stat)=\aleph_2$ implies $\continuum\leq\aleph_2$. 
  (\CH\ is also possible under $\lsn(\calL^{\aleph_0}_\stat)=\aleph_2$).\smallskip

  \wassert{7} 
  If $\kappa\in\LSS(\calL^{\aleph_0}_\stat)$ and
  $\kappa>\aleph_2$ then $\kappa>\continuum$. \smallskip

  \wassert{8} 
  $\lsn(\calL^{\aleph_0}_\stat)=\continuum$ implies
  $\continuum=\aleph_2$. \smallskip
  
  \wassert{9} If $\kappa$ is $\sigma$-closed-gen.\ supercompact then
  $\lsn(\calL^{\aleph_0,\II}_\stat)\leq\kappa$.\smallskip

  \wassert{10} Assuming the consistency of \ZFC\ $+$ ``there is a supercompact cardinal'', 
  the following assertion is consistent with \ZFC\ for any natural number $n\geq2$:
  \begin{xitemize}
  \xitem[x-stat-0] $\lsn(\calL^{\aleph_0,\II}_\stat)=\aleph_n$. 
  \end{xitemize}
\end{Prop}
\prf \assertof{1}: is easy to prove\ifextended\else\ (for the details of a proof see [0])\fi.
\ifextended{\extendedcolor
Let 
$\calC:=\setof{\kappa\in\Card}{\kappa\mbox{ is }{\cdot}^{\aleph_0}\mbox{-closed\ }}$. 

$\calC\subseteq\LSS(\calL^{\aleph_0,\II})$:\quad Suppose that $\kappa$ 
is $\cdot^{\aleph_0}$-closed. Let $\gmA$ be an arbitrary structure of countable signature 
and $S\in[\,\gmA\,]^{\LT\kappa}$. Let $\theta$ be a sufficiently large regular cardinal
($\gmA\in\calH(\theta)$ in particular).

Let $M\prec\calH(\theta)$ be \st\ $\gmA\in M$, $S\subseteq M$, $[M]^{\aleph_0}\subseteq M$ 
and $\cardof{M}<\kappa$ (this is possible since $\kappa$ is $\cdot^{\aleph_0}$-closed).
Then $\gmA\restr (\ulsetof{\gmA}\cap M)$ contains $S$ as a subset of its underlying set
$\ulsetof{\gmA}\cap M$, $\cardof{\ulsetof{\gmA}\cap M}<\kappa$, and
$\gmA\restr (\ulsetof{\gmA}\cap M)\prec_{\calL^{\aleph_0,\II}}\gmA$.
\smallskip

$\LSS(\calL^{\aleph_0,\II})\subseteq\calC$:\quad Suppose that $\kappa$ is not in $\calC$ and 
let $\mu<\kappa$ be \st\ $\mu^{\aleph_0}\geq\kappa$.

Let $\gmA:=\pairof{\calH(\mu^+),\in}$. Note that
\begin{xitemize}
\item[] 
  $\gmA\models_{\calL^{\aleph_0,\II}}\forall X\exists_x\forall_y\,(y\varin X\ \leftrightarrow\ y\in x)$.
\end{xitemize}

If $\gmB\prec_{\calL^{\aleph_0,\II}}\gmA$, then
\begin{xitemize}
\item[] 
  $\gmB\models_{\calL^{\aleph_0,\II}}\forall X\exists_x\forall_y\,(y\varin X\ \leftrightarrow\ y\in x)$
\end{xitemize}
by elementarity, and hence $[\,\ulsetof{\gmB}\,]^{\aleph_0}\subseteq\ulsetof{\gmB}$. Thus, 
if $\mu\subseteq\ulsetof{\gmB}$ then $[\mu]^{\aleph_0}\subseteq\ulsetof{\gmB}$, and hence 
we have $\kappa\leq\mu^{\aleph_0}\leq\ccardof{\gmB}$. This shows that
$\kappa\not\in\LSS(\calL^{\aleph_0,\II})$. 
}\fi 
\smallskip

\assertof{2}: follows from \assertof{1}.\smallskip

\assertof{3}: If $\square_\kappa$ holds then there is a non-reflecting stationary set
$S\subseteq E^{\kappa^+}_\omega$. 

$S^*:=\setof{s\in[\kappa^+]^{\aleph_0}}{\sup(\setof{\alpha+1}{\alpha\in s})\in S}$ is a non-reflecting 
stationary subset of $[\kappa^+]^{\aleph_0}$. $S^*$ can be recast to a counter-example of 
$\lsn(\calL^{\aleph_0}_{\stat})\leq\kappa^+$. 
\smallskip

\assertof{4}: follows from \assertof{3}.\smallskip

\assertof{5}: This can be proved by a modification of arguments in \cite{sfetal-I}. 
The second assertion holds since $\MA^{+\omega_1}(\sigma\mbox{-closed})$ $+$
$\continuum=\aleph_2$ is consistent (e.g.\ modulo a supercompact). 
\smallskip

\assertof{6}: $\lsn(\calL^{\aleph_0}_\stat)=\aleph_2$ implies the principle \RP\ in 
\cite{millennium-book}. \RP\ implies $\continuum\leq\aleph_2$  (Todor\v{c}evi\'c, see 
Theorem 37.17 in \cite{millennium-book}).

Suppose $\uniV$ satisfies \CH\ and there is a 
supercompact cardinal $\kappa$ then collapsing $\kappa$ to $\aleph_2$, we obtain a model of
$\lsn(\calL^{\aleph_0}_\stat)=\aleph_2$ (see \assertof{9}) and \CH. 
\smallskip

\assertof{7}: If $\kappa\in\lsn(\calL^{\aleph_0}_\stat)$ then
$\SDLS^-(\calL^{\aleph_0}_\stat,\,\LT\kappa)$ holds (see \cite{sfetal-I} for definition). 
Thus,  by Proposition 2.1 in \cite{sfetal-II}, if $\kappa>\aleph_2$ then
$\kappa>\continuum$ holds. 
\smallskip

\assertof{8}: This follows from Corollary 2.3 in \cite{sfetal-II}. 
\smallskip

\assertof{9}: See \cite{sfetal-I} (for the detail of a direct proof, see [0]). 

\ifextended{\extendedcolor
Suppose that $\gmA$ is a structure with $\ccardof{\gmA}\geq\kappa$ and
$S\in[\,\ulsetof{\gmA}\,]^{\LT\kappa}$. Let $\poP$ be a $\sigma$-closed \po\ with
a $(\uniV,\poP)$-generic $\genG$ \st\ there are $j,M\subseteq\uniV[\genG]$ with
\begin{xitemize}
\xitemA[x-LS-25] 
  $\Elembed{j}{\uniV}{M}{\kappa}$,
\xitemA[x-LS-26] 
  $j(\kappa)>\ccardof{\gmA}$ and
\xitemA[x-LS-27] 
  $j\imageof \mu\in M$ where $\mu:=\ccardof{\gmA}^{\aleph_0}$.
\end{xitemize}
\Wolog, we may assume that \ixitem[x-LS-27-a] $\ulsetof{\gmA}=\ccardof{\gmA}$. 

Let
$\gmB:=j(\gmA)\restr j\imageof{\ulsetof{\gmA}}$. By \xitemof{x-LS-27} and \xitemof{x-LS-27-a}, $\gmB$,
$j\restr\ulsetof{\gmA}\in M$. Thus, $\gmA\in M$ and we have
\begin{xitemize}
\xitemA[x-LS-27-0] 
  $M\modelof{\isom{j\restr\ulsetof{\gmA}}{\gmA}{\gmB}}$.
\end{xitemize}
Since
$([\,\ulsetof{\gmA}\,]^{\aleph_0})^\uniV=([\,\ulsetof{\gmA}\,]^{\aleph_0})^{\uniV[\genG]}$ 
by $\sigma$-closedness of $\poP$ and
$([\,\ulsetof{\gmA}\,]^{\aleph_0})^\uniV\subseteq([\,\ulsetof{\gmA}\,]^{\aleph_0})^M$, we 
have
\begin{xitemize}
\xitemA[x-LS-28] $([\,\ulsetof{\gmA}\,]^{\aleph_0})^\uniV=([\,\ulsetof{\gmA}\,]^{\aleph_0})^M$.
\end{xitemize}
\begin{Claim}
  \Label{Cl-LS-8} For any $\calL^{\aleph_0,\II}_\stat$-formula
  $\varphi(x_0\ctentenc X_0\ctenten)$,  $a_0\ctenten\in\ulsetof{\gmA}$ and
  $A_0\ctenten\in[\,\ulsetof{\gmA}\,]^{\aleph_0}$, 
  \begin{xitemize}
  \xitemA[x-LS-29] 
    $\uniV\modelof{\gmA\models\varphi(a_0\ctentenc A_0\ctenten)}$\ \ $\Leftrightarrow$\ \
    $M\modelof{\gmA\models\varphi(a_0\ctentenc A_0\ctenten)}$. 
  \end{xitemize}
\end{Claim}
\prfofClaim
We prove \xitemof{x-LS-29} by induction on $\varphi$. The most crucial step is when 
$\varphi$ is of the form $stat\, X\,\psi(x_0\ctentenc X, X_0\ctenten)$:

Suppose that $a_0\ctenten\in\ulsetof{\gmA}$ and
$A_0\ctenten\in[\,\ulsetof{\gmA}\,]^{\aleph_0}$.

First assume that $\uniV\modelof{\gmA\models stat\, X\,\psi(a_0\ctentenc X, A_0\ctenten)}$. 
This is equivalent to 
\begin{xitemize}
\item[] $\uniV\modelof{\setof{A\in[\,\ulsetof{\gmA}\,]^{\aleph_0}}{
    \gmA\models\psi(a_0\ctentenc A, A_0\ctenten)}\mbox{ is stationary}}$.
\end{xitemize}
We have 
\begin{xitemize}
\item[] 
  $(\setof{A\in[\,\ulsetof{\gmA}\,]^{\aleph_0}}{
  \gmA\models\psi(a_0\ctentenc A, A_0\ctenten)})^\uniV\\[\jot]
  \qquad=(\setof{A\in[\,\ulsetof{\gmA}\,]^{\aleph_0}}{
    \gmA\models\psi(a_0\ctentenc A, A_0\ctenten)})^M$ 
\end{xitemize}
by \xitemof{x-LS-28} and induction hypothesis. Since $\poP$ is proper, it follows that
$M\modelof{\setof{A\in[\,\ulsetof{\gmA}\,]^{\aleph_0}}{
    \gmA\models\psi(a_0\ctentenc A, A_0\ctenten)}\mbox{ is stationary}}$. Thus 
$M\modelof{\gmA\models stat\, X\,\psi(a_0\ctentenc X,A_0\ctenten)}$.

Now assume that $\uniV\not\modelof{\gmA\models stat\, X\,\psi(a_0\ctentenc X,A_0\ctenten)}$.
This means that there is a club $C\subseteq[\,\ulsetof{\gmA}\,]^{\aleph_0}$ \st\
$\uniV\modelof{\mbox{for all }A\in C,\,\gmA\not\models\psi(a_0\ctentenc A,A_0\ctenten)}$. 
Now $C\in M$ by \xitemof{x-LS-27}. By induction hypothesis 
\begin{xitemize}
\item[]
  $M\modelof{\mbox{for all }A\in C,\,\gmA\not\models\psi(a_0\ctentenc A,A_0\ctenten)}$. 
\end{xitemize}
By $\sigma$-closedness of $\poP$, $C$ is club in $\uniV[\genG]$ and hence also in $M$. 

Thus, it follows that $M\not\modelof{\gmA\models stat\, X\,\psi(a_0\ctentenc X,A_0\ctenten)}$. 
\qedofClaim
\qedskip

For an $\calL^{\aleph_0,\II}_\stat$-formula $\varphi(x_0\ctentenc X_0\ctenten)$,
$a_0\ctenten\in\ulsetof{\gmA}$ and $A_0\ctenten\in[\,\ulsetof{\gmA}\,]^{\aleph_0}$
\begin{xitemize}
\item[] $M\modelof{j(\gmA)\models\varphi(j(a_0)\ctentenc j(A_0)\ctenten)}$
  \ \ $\obecause{\Leftrightarrow}{}{by elementarity of $j$}$\ \ 
  $\uniV\modelof{\gmA\models\varphi(a_0\ctentenc A_0\ctenten)}$\\[\jot]
  $\ubecause{\Leftrightarrow}{}{by \Claimof{Cl-LS-8}}$\ \ 
  $M\modelof{\gmA\models\varphi(a_0\ctentenc A_0\ctenten)}$\\
  $\ubecause{\Leftrightarrow}{}{by \xitemof{x-LS-27-0}}\ \ 
  M\modelof{\gmB\models\varphi(j(a_0)\ctentenc j(A_0)\ctenten)}$.
\end{xitemize}

Since $M\models\ccardof{\gmB}<j(\kappa),\,j(S)\ (=j\imageof S)\ \subseteq \ulsetof{\gmB}$, we obtain
\begin{xitemize}
\item[] $M\modelof{
  \begin{array}[t]{@{}l}
    \mbox{there is an }\calL^{\aleph_0,\II}_\stat\mbox{-elementary substructure }X
    \mbox{ of }j(\gmA)
    \mbox{ of size }\LT j(\kappa)\\
    \mbox{with }j(S)\subseteq X}.
  \end{array}
  $
\end{xitemize}
By elementarity, it follows that 
\begin{xitemize}
\item[] $\uniV\modelof{
  \begin{array}[t]{@{}l}
    \mbox{there is an }\calL^{\aleph_0,\II}_\stat\mbox{-elementary substructure }X
    \mbox{ of \,}\gmA
    \mbox{ of size }\LT\kappa\mbox{ with }\\S\subseteq X}. 
  \end{array}
  $
\end{xitemize}

This shows that $\kappa\in\LSS(\calL^{\aleph,\II}_\stat)$. 
}\fi 
\smallskip

\assertof{10}: By \assertof{3} and \assertof{5}. Note that, if $\kappa$ is a supercompact 
cardinal, and $\poP=\Col(\aleph_{k},\kappa)$ for $k>0$ (in the notation of Kanamori \cite{higher-inf}) 
then, in $\uniV[\genG]$ for a $(\uniV,\poP)$-generic $\genG$, 
$\kappa=\aleph_{k+1}$ and $\kappa$ is $\sigma$-closed-gen.\ supercompact. 
\qedofProp\qedskip

In contrast to the possible smallness of $\lsn(\calL^{\aleph_0,\II})$ shown in 
\Propof{P-stat-0},\,\assertof{2} (e.g. it can be $\aleph_2$ under \CH),
$\wcn(\calL^{\aleph_0,\II})$ is fairly a large cardinal. This is because we can apply \Propof{P-comp-4-0} 
to $\calL^{\aleph_0,\II}$. 




For a logic $\calL$, the Hanf number $\darkred\hn(\calL)$ of $\calL$ is defined by
\begin{xitemize}
\xitem[x-stat-1]
  ${\darkred\hn(\calL)}:=\min(\setof{\kappa}{
  \begin{array}[t]{@{}l}
    \mbox{for any }\calL\mbox{-sentence }\varphi
    \mbox{ if }\varphi\mbox{ has a model of size at least }\\
    \kappa\mbox{ then }\varphi\mbox
         { has arbitrarily large model}}).
  \end{array}$
\end{xitemize}

Note that all proper logics have Hanf number assuming the properness include the 
properties that, there is a cardinal $\kappa_\calL$ \st\ each $\calL$-formula is in a signature of 
size $\LT\kappa_\calL$ and there are only set-many $\calL$-formulas in each signature. 
The following should be well-known:
\begin{Lemma}\Label{P-stat-2}
  For a proper logic $\calL$, we have $\hn(\calL)\leq\cn(\calL)$. \ifextended\else\qed\fi
\end{Lemma}
\ifextended{\extendedcolor
\prf Suppose that $\kappa<\hn(\calL)$ then there is an $\calL$-sentence $\varphi$ \st\ 
$\varphi$ has a model $\gmA$ with $\kappa\leq\ccardof{\gmA}<\hn(\calL)$ \st\ any model of 
$\varphi$ has size $\LT\hn(\calL)$.

Let
$T:=\ssetof{\varphi}\cup\setof{\symbc_\alpha\not\equiv\symbc_\beta}{\alpha<\beta<\hn(\calL)}$. 
Then $T$ does not have any model but, for any
$T_0\in[T]^{<\kappa^+}$, $\gmA$ can be expanded to a model of $T_0$. This shows that
$\cn(\calL)>\kappa^+$. (This proof actually shows that $\hn(\calL)<\cn(\calL)$ if 
$\hn(\calL)$ is a successor. But this fact is irrelevant for the application in the 
following Lemma). 
\qedofLemma\qedskip}\fi 

In the following \Lemmaof{P-stat-3}, $Q_1$ denotes the (first-order) quantifier ``there exists uncountably 
many'' and $\calL(Q_1)$ the logic obtained by adding the quantifier $Q_1$ to the first-order logic. 

The next lemma follows immediately from \Lemmaabove. Note the quantifier $Q_1$ is 
interpretable in $\calL^{\aleph_0}_\stat$. \ifextended\\{\extendedcolor
{[\,} ${Q_1}_{\,x\,}\varphi$\ \ $\mapsto$\ \
  $stat\, X\exists_x(x\notvarin X\,\land\,\varphi^*)$ where $\varphi^*$ is the 
$\calL^{\aleph_0}_{\stat}$-interpretation of $\varphi$. {\,]}
}\fi 

\begin{Lemma}
  \Label{P-stat-3}\wassertof{1}\
  $\wcn(\calL^{\aleph_0}_\stat)\geq\wcn(\calL(Q_1))\geq\aleph_\omega$.\smallskip

  \wassert{2} {\rm(Shelah \cite{ShelahII})} For a proper Logic $\calL$, $\hn(\calL)$ is a 
  strong limit. If $\hn(\calL)>\aleph_0$, then $\cn(\calL)\geq\hn(\calL)\geq\beth_\omega$.\smallskip

  \wassert{3} {\rm(Shelah \cite{shelah-kaufmann})}\quad $\hn(\calL^{\aleph_0}_{\stat})>\beth_\omega$. 
  In particular $\cn(\calL^{\aleph_0}_\stat)>\beth_\omega$, while we have
  $\cn(\calL^{\aleph_0}_\stat)\geq\cn(\calL(Q_1))\geq\hn(\calL(Q_1))=\beth_\omega$. \qed
\end{Lemma}
\memo{Is $\wcn(\calL^{\aleph_0}_\stat)<\continuum$ consistent?\\
  How about $\cn(\calL^{\aleph_0}_\stat)<\continuum$?\\
  $\wcn(...)<\cn(\calL^{\aleph_0}_\stat)$?\\
  What about the Cohen model over $\uniL$?}



For the \Propof{P-stat-4} below, we use Jensen's global square $\globalC$ (for the 
existence of the class $\globalC$ as below, see the proof of Theorem 5.1 in \cite{jensen}):

Let 
\begin{xitemize}
\item[] ${\darkred\SINGL}:=\setof{\alpha}{\uniL\modelof{\alpha\mbox{ is a singular limit ordinal}}}$.
\end{xitemize}

A global square is a class function $\globalC$ (i.e.\ the corresponding $\Lin$-formula defining 
the class function) \st\ the following are provable from $\ZFC^-$:

\begin{xitemize}
\item[\wassertof{I}\ \ \ ]\!\! $\globalC$ is a class function with $\dom(\globalC)=\SINGL$.
\item[\wassertof{II}\ \ \ ]\!\! For each $\alpha\in\SINGL$, $\globalC(\alpha)\in L$ is a club subset of 
  $\alpha$ with $\otp(\globalC(\alpha))<\alpha$.
\item[\wassertof{III}\ \ \ ]\!\! {\rm(Coherence)}\quad If $\alpha\in\SINGL$, and
  $\beta\in\Lim(\globalC(\alpha))$, then $\beta\in\SINGL$ and $\globalC(\beta)=\globalC(\alpha)\cap\beta$. 
\item[\wassertof{IV}\ \ \ ]\!\! For any transitive set model $W$ of $\ZFC^-$, we have
  $\globalC^W=\globalC\restr(\SINGL)^W$. 
\end{xitemize}

\begin{Lemma}\Label{P-stat-3-0}{\rm (see e.g.\ Kunen \cite{kunen-book}, Lemma 4.11)}
  Assume $\uniV = \uniL$. Then $L_\kappa=\calH(\kappa)$ whenever $\kappa>\omega$ and $\kappa$ is regular.\qed
\end{Lemma}

\begin{Prop}
  \Label{P-stat-4} 
  Assume that $\kappa=\wcn(\calL^{\aleph_0,\II}_{\stat, \kappa,\omega})$ and  
  suppose that $\lambda>\kappa$ is a regular cardinal and $M\prec L_\lambda$ is 
  \st\ $\kappa+1\subseteq M$ and $\cardof{M}=\kappa$. Then there are 
  $\overline{\lambda}>\kappa$ and $j$ \st\ 
  \begin{xitemize}
  \xitem[x-stat-1-a] $\overline{\lambda}$ is a regular cardinal in $\uniL$ of uncountable 
    cofinality (in $\uniV$), and
  \xitem[x-stat-1-a-0] $\EElembed{j}{M}{L_{\overline{\lambda}}}{\kappa}$\,. 
  \end{xitemize}
\end{Prop}
\prf If $0^\#$ exists, then, since the regular cardinals $\kappa$, $\lambda$ belong to the 
indiscernibles associated with $0^\#$, we can easily find $j$ as above for any regular
$\overline{\lambda}\geq\kappa^+$.

In the following, we assume that $0^\#$ does not exist. Observe first that it is enough to 
prove the following variation of the statement of the present proposition:
\begin{xitemize}
\xitem[x-stat-2-0] For any $M\prec L_\lambda$ with $\kappa+1\subseteq M$ and
  $\cardof{M}=\kappa$, there is $M^*\prec L_\lambda$ \st\ $M\subseteq M^*$ and
  $\cardof{M^*}=\kappa$, with $\overline{\lambda}>\kappa$ and $j^*$ \st\ 
  \begin{xitemize}
  \item[\xitemof{x-stat-1-a} \ \ \ ]\!\!\!$\overline{\lambda}$ is a regular cardinal, and
  \item[\xitemdof{x-stat-1-a-0}\ \ \ ]\!\!\!$\EElembed{j^*}{M^*}{L_{\overline{\lambda}}}{\kappa}$\,. 
  \end{xitemize}
\end{xitemize}

This is because, if $M^*$, $\overline{\lambda}$, $j^*$ are as in \xitemof{x-stat-2-0} for 
given $M$, then $M$, $\overline{\lambda}$ and $j:=j^*\restr M$ satisfy \xitemof{x-stat-1-a} 
and \xitemof{x-stat-1-a-0}. 

Let $E:=(E^\lambda_\omega\setminus\omega_2)^\uniV$. For any $\alpha\in E$, 
we have $(\cf(\alpha))^\uniL<\omega_2$ by Covering Lemma. Thus, $\globalC(\alpha)$ exists for 
such $\alpha$ and $\otp(\globalC(\alpha))<\alpha$ by \assertof{II}. By Fodor's Lemma, there is 
$\nu<\lambda$ \st\ 
\begin{xitemize}
\xitem[x-stat-2-1] 
  $E_0:=\setof{\alpha\in E}{\otp(\globalC(\alpha))=\nu}$ is stationary. 
\end{xitemize}
Since $\cf(\nu)=\omega$ (in $\uniV$), this statement is equivalent to 
\begin{xitemize}
\xitem[x-stat-2-2]
  $E_1:=\setof{a\in[\lambda]^{\aleph_0}}{a\mbox{ does not have the maximal element, }
  \sup(a)\in E_0}$ is stationary subset of $[\lambda]^{\aleph_0}$.
\end{xitemize}

Now suppose that $M\prec L_\lambda$ is \st\ $\kappa+1\subseteq M$ and $\cardof{M}=\kappa$. 
Let $M^*\prec L_\lambda$ be \st\ $M\cup\ssetof{\nu}\subseteq M^*$ 
and $\cardof{M^*}=\kappa$.

Let
\begin{xitemize}
\xitem[x-stat-2-3] $T:={}
  \begin{array}[t]{@{}l}
    \obecause{\setof{\varphi(\symbc_{a_0}\ctenten)}{
      \begin{array}[t]{@{}l}
        a_0\ctenten\in M^*,\varphi(x_0\ctenten)\mbox{ is an }
        \Lstatinf\mbox{-formula with}\\L_\lambda\models\varphi(a_0\ctenten)}
    \end{array}}{}{\normalsize$(\ast)$}\\[\jot]
    \cup\ \setof{\symbc_\alpha<\symb{d}}{\alpha\in\kappa}\\
    \cup\ \ssetof{\symb{d}<\symbc_\kappa}\\
    \cup\ \setof{\forall_x\,((x\mbox{ is an ordinal}\land x<\symbc_\alpha) 
      \rightarrow \llor_{\beta<\alpha}x\equiv\symbc_\beta)}{\alpha<\kappa}.
    \ifarxived 
    \smash{\raisebox{3.24ex}{$\left.\rule[0cm]{0cm}{4.68ex}\right\}$}\raisebox{3.42ex}{$(\ast\ast)$}}
    \else
    \smash{\raisebox{1.8zh}{$\left.\rule[0cm]{0cm}{2.6zh}\right\}$}\raisebox{1.9zh}{$(\ast\ast)$}}
    \fi
  \end{array}
  $
\end{xitemize}

Apparently $T$ is of a signature of size $\kappa$. 
As before we can also show that $T$ is $\LT\kappa$-satisfiable. 
Hence, by
$\kappa=\wcn(\Lstatinf)$, $T$ has a model $\gmA$.

Since $\symb{\in}^\gmA$ is well-founded (which is declared in the $(\ast)$ part of $T$), 
we can take the Mostowski collapse $\gmA^*$ of $\gmA$. The underlying set of $\gmA^*$ is of 
the form $L_{\overline{\lambda}}$ (also because this is formulated in the $(\ast)$ part of 
$T$).

Letting $\mapping{j^*}{M^*}{L_{\overline{\lambda}}}$; $a\mapsto(\symbc_a)^{\gmA^*}$, we have
\begin{xitemize}
\xitem[x-stat-2-4] $\EElembed{j^*}{M^*}{L_{\overline{\lambda}}}{\kappa}$
\end{xitemize}
(the elementarity follows from $(\ast)$ part of $T$, and  $\crit(j)=\kappa$ from $(\ast\ast)$ of $T$), and 
\begin{xitemize}
\xitem[x-stat-2-5] 
  $\pairof{L_\lambda, a}_{a\in M^*}\equiv_{\Lstatinf}
  \pairof{L_{\overline{\lambda}}, j^*(a)}_{a\in M^*}$ 
\end{xitemize}
by the $(\ast)$ part of $T$.

Note that $\cf(\overline{\lambda})>\omega$ by 
\xitemof{x-stat-2-5}. 
Thus, we are done by showing that $\overline{\lambda}$ is regular in $\uniL$. Suppose not, by way of 
contradiction. 
Then $\globalC(\overline{\lambda})$ exists. 
$C:=\setof{\alpha<\overline{\lambda}}{\alpha\mbox{ is a limit of }\globalC(\overline{\lambda})}$ 
is a club subset of $\overline{\lambda}$. 
$E_2:=\setof{\alpha<\overline{\lambda}}{\otp(\globalC(\alpha))=j^*(\nu)}$ is stationary by 
\xitemof{x-stat-2-2}, the elementarity \xitemof{x-stat-2-5}) and \wassertof{IV}.

Let $\xi_0$, $\xi_1\in C\cap E_2$ be with $\xi_0<\xi_1$. Then $\globalC(\xi_0)=\globalC(\xi_1)\cap\xi_0$ 
by \assertof{III}. This is a contradiction to $\xi_0$, $\xi_1\in E_2$. 
\qedofProp
\qedskip

Under $\uniV=\uniL$, the 
condition $\kappa=\wcn(\Lstatinf)$ becomes equivalent to
$\kappa=\wcn(\calL^\II_{\kappa,\omega})$: 

\begin{Thm}
  \phantomsection
  \Label{P-stat-7} Assume  $\uniV=\uniL$. Then  
  $\kappa=\wcn(\Lstatinf)$ holds if and only if $\kappa$ is weakly extendible. 
\end{Thm}
\memo{\texttt{weak\char`_stat\char`_compact 4.pdf}}

\Thmof{P-stat-7} follows from \Propof{P-stat-4} and the following characterization of weak 
extendibility: 
\begin{Lemma}
  \Label{P-stat-8} For a cardinal $\kappa$, \tfae:\smallskip

  \wassert{a} $\kappa$ is weakly extendible.\smallskip

  \wassert{b} $2^{\LT\kappa}=\kappa$ holds, and, for any regular $\lambda>\kappa$ and any 
$M\prec\calH(\lambda)$ with $\kappa+1\subseteq M$ and $\cardof{M}=\kappa$, there are regular 
$\overline{\lambda}>\kappa$ and $\EElembed{j}{M}{\calH(\overline{\lambda})}{\kappa}$. 
\end{Lemma}
\prf \assertof{a} $\Rightarrow$ \assertof{b}: Assume that $\kappa$ is weakly extendible. 
Then we have $2^{\LT\kappa}=\kappa$ by definition. For a regular $\lambda>\kappa$ and 
$M\prec\calH(\lambda)$ with $\kappa+1\subseteq M$ and $\cardof{M}=\kappa$, let
$\theta:=\lambda+\omega$. Then $M$, $\calH(\lambda)\in V_\theta$. Let $M^*$ be \st\ 
\begin{xitemize}
\xitem[x-stat-9] $M^*\prec V_\theta$, 
\xitem[x-stat-10] $M$, $\calH(\lambda)\in M^*$, $M\subseteq M^*$ and
\xitem[x-stat-11] $\cardof{M^*}=\kappa$. 
\end{xitemize}
Since $\kappa$ is weakly extendible, there are $\overline{\theta}$ and $\overline{j}$ with
\ixitem[x-stat-12] $\EElembed{\overline{j}}{M^*}{V_{\overline{\theta}}}{\kappa}$. Note 
that $\overline{\theta}$ is a limit ordinal by the choice of $\theta$ and by the elementarity 
\xitemof{x-stat-9} and \xitemof{x-stat-12}. 

We have
$V_{\overline{\theta}}\modelof{\overline{j}(\lambda)\mbox{ is a regular cardinal}}$  by 
elementarity \xitemof{x-stat-12}, and hence $\overline{j}(\lambda)$ is really a regular 
cardinal. By the elementarity, 
we also 
have $V_{\overline{\theta}}\modelof{\overline{j}(\calH(\lambda))=\calH(\overline{j}(\lambda))}$. Since
$\calH(\overline{j}(\lambda))\in V_{\overline{\theta}}$, it follows that
$\overline{j}(\calH(\lambda))=\calH(\overline{j}(\lambda))$.

Thus, letting $j:=\overline{j}\restr M$ and $\overline{\lambda}:=\overline{j}(\lambda)$,  we 
obtain $\EElembed{j}{M}{\calH(\overline{\lambda})}{\kappa}$. \smallskip

\assertof{b} $\Rightarrow$ \assertof{a}: Assume that \assertof{b} holds. We want to show 
that \xitemof{x-w-ext-1} holds. Suppose $\theta>\kappa$ and $M\prec V_\theta$ is \st\
$\kappa+1\subseteq M$ and $\cardof{M}=\kappa$.

Let $\lambda>\theta$ be a regular cardinal \st\ $V_\theta\in\calH(\lambda)$, and let 
$M^*\prec\calH(\lambda)$ be \st\ $M$, $\theta$, $\lambda\in M^*$, $M\subseteq M^*$ and
$\cardof{M^*}=\kappa$. 

By \assertof{b}, there are regular $\overline{\lambda}>\kappa$ and $j$ \st\
$\EElembed{\overline{j}}{M^*}{\calH(\overline{\lambda})}{\kappa}$. Then, by letting
$j:=\overline{j}\restr M$ and $\overline{\theta}:=\overline{j}(\theta)$,  
we obtain $\EElembed{j}{M}{V_{\overline{\theta}}}{\kappa}$. 
\qedofLemma\qedskip

\noindent
{\bf Proof of \bfThmof{P-stat-7}:} If $\kappa$ is weakly extendible then 
$\kappa=\wcn(\calL^\HO_{\kappa,\kappa})$ holds (\Thmof{P-w-ext-0}). By 
\Lemmaof{P-intro-8},\,\assertof{3}, it follows that $\kappa=\wcn(\Lstatinf)$. 

Now assume that $\kappa=\wcn(\Lstatinf)$. We want to show that \assertof{b} in  
\Lemmaof{P-stat-8} holds. 

We have $2^{\LT\kappa}=\kappa$ since $\kappa=\wcn(\calL_{\kappa,\omega})$ and hence 
$\kappa$ is weakly compact.

Suppose that $\lambda>\kappa$ is regular and $M\prec\calH(\lambda)$ is \st\
$\kappa+1\subseteq M$ and $\cardof{M}=\kappa$. By \Lemmaof{P-stat-3-0}, we have
$\calH(\lambda)=L_\lambda$. Thus, by \Propof{P-stat-4}, there are regular 
$\overline{\lambda}>\kappa$ and $j$ with $\EElembed{j}{M}{L_{\overline{\lambda}}}{\kappa}$.
Again by \Lemmaof{P-stat-3-0}, we have
$\calH(\overline{\lambda})=L_{\overline{\lambda}}$. This shows that \assertof{b} in 
\Lemmaof{P-stat-8} holds. 
\qedof{\Thmof{P-stat-7}}\qedskip

For a cardinal $\kappa$ with $\kappa=\wcn(\Lstatinf)$, we do not know if the same equation 
holds in $\uniL$. However we can show that quite strong large cardinal properties hold 
around $\kappa$ in $\uniL$: 

\begin{Thm}
  \Label{P-stat-6} Assume that $\kappa=\wcn(\Lstatinf)$. 
  Then we have:\smallskip

  \wassert{1} $\uniL\modelof{{}
    \begin{array}[t]{@{}l}
      \kappa\mbox{ is weakly compact, and there are stationarily many weakly compact}\\
      \mbox{cardinals below }\kappa}.
    \end{array}$
  \smallskip

  \wassert{2} $\uniL\modelof{{}
    \begin{array}[t]{@{}l}
      \mbox{there is a weakly compact cardinal }\GT\kappa\mbox{ which is a limit of weakly}\\
      \mbox{compact cardinals}}.
    \end{array}$
\end{Thm}
\prf \assertof{1}: As already noted in the proof of \Thmof{P-stat-7}, $\kappa$ is weakly 
compact. It follows that $\kappa$ is also weakly compact in $L$ (see e.g.\ Theorem 17.22 in 
\cite{millennium-book}).

Suppose that $D\in L$ is a closed unbounded subset of $\kappa$. We have to show that there 
is $\mu\in D$ \st\ $\uniL\modelof{\mu\mbox{ is weakly compact}}$. 

Let $\lambda>\kappa$ be a regular cardinal and let \ixitem[x-stat-13] $M\prec L_\lambda$ be 
\st\ $\kappa+1\subseteq M$, $D\in M$, and $\cardof{M}=\kappa$. By \Propof{P-stat-4}, there 
are $\overline{\lambda}>\kappa$ regular in $\uniL$, and $j$ \st\
\ixitem[x-stat-14] $\EElembed{j}{M}{L_{\overline{\lambda}}}{\kappa}$. Since
$\overline{\lambda}$ is regular in $\uniL$, we have
$L_{\overline{\lambda}}\modelof{\kappa\mbox{ is weakly compact}}$ (this can be seen e.g.\ in the 
characterization of weak compactness in terms of the tree property).
By the elementarity 
\xitemof{x-stat-13} and \xitemof{x-stat-14}, and since $D=j(D)\cap\kappa$ is unbounded in
$\kappa$, we have $L_{\overline{\lambda}}\modelof{\kappa\in j(D)}$. Thus 
\begin{xitemize}
\item[] 
  $L_{\overline{\lambda}}\modelof{\mbox{there is weakly compact }\mu\in j(D)}$.
\end{xitemize}
Hence $L_\lambda\modelof{\mbox{there is weakly compact }\mu\in D}$ by the elementarity 
\xitemof{x-stat-13} and \xitemof{x-stat-14}. Since $\lambda$ is regular, it follows that
$\uniL\modelof{\mbox{there is weakly compact }\mu\in D}$.\smallskip

\assertof{2}: Let $\lambda>\kappa$ be regular and \ixitem[x-stat-15] $M\prec L_\lambda$ be 
\st\ $\kappa+1\subseteq M$ and $\cardof{M}=\kappa$. By \Propof{P-stat-4}, there are regular
$\overline{\lambda}>\kappa$ and $j$ with \ixitem[x-stat-16]
$\EElembed{j}{M}{L_{\overline{\lambda}}}{\kappa}$. 
By \assertof{1} and since $\lambda>\kappa$ is regular, we have
\begin{xitemize}
\item[] 
  $L_\lambda\modelof{{}
  \begin{array}[t]{@{}l}
    \kappa\mbox{ is weakly compact and is stationary limit of weakly compact}\\
    \mbox{cardinals}.}
  \end{array}$
\end{xitemize}

By the elementarity \xitemof{x-stat-15} and \xitemof{x-stat-16}, it follows that 
\begin{xitemize}
\item[] 
  $L_{\overline{\lambda}}\modelof{{}
  \begin{array}[t]{@{}l}
    j(\kappa)\mbox{ is weakly compact and is stationary limit of weakly compact}\\
    \mbox{cardinals}.}
  \end{array}$
\end{xitemize}
Since $\overline{\lambda}$ is regular in $\uniL$, it follows that 
\begin{xitemize}
\item[] 
  $\uniL\modelof{{}
  \begin{array}[t]{@{}l}
    j(\kappa)\mbox{ is weakly compact and is stationary limit of weakly compact}\\
    \mbox{cardinals}.}
  \end{array}$\ifarxived\vspace{-2.16ex}\else\vspace{-1.2zh}\fi
\end{xitemize}

\qedofThm
\qedskip

There are still many open problems concerning the compactness of weak second-order logic 
and stationary logic including the one we mentioned before \Thmof{P-stat-6}. We neither 
know the answer to the following:

\begin{Problem}
  Is ``\/$\wcn(\calL^{\aleph_0,\II}_\stat)<\wcn(\calL^\II)$'' consistent? 
\end{Problem}

\noindent
{\bf Acknowledgments: } The authors would like to thank Philipp L\"ucke for pointing out a 
serious flaw in an early draft of the paper.

\end{document}